\documentclass[12pt]{amsart}
\usepackage{amssymb,amsmath,mathrsfs}
\usepackage[ps2pdf,colorlinks=true,urlcolor=blue,
citecolor=red,linkcolor=blue,linktocpage,pdfpagelabels,bookmarksnumbered,bookmarksopen]{hyperref}
\usepackage{epsfig,graphicx,color}
\usepackage[english]{babel}

\usepackage[left=2.6cm,right=2.6cm,top=2.6cm,bottom=2.6cm]{geometry}

\newcommand{\N}{{\mathbb N}}

\newcommand{\C}{{\mathbb C}}
\newcommand{\R}{{\mathbb R}}

\newcommand{\rn}{{\mathbb{R}^N}}
\newcommand{\be}{\begin{equation}}
\newcommand{\ee}{\end{equation}}

\renewcommand{\H}{{\mathbb H}}

\newcommand{\iu}{{\rm i}}
\newcommand{\teta }{\theta }
\newcommand{\vep}{\varepsilon}
\newcommand{\eps}{\varepsilon}

\newcommand{\E}{{\mathcal E}}

\newcommand{\cala}{{\mathcal A}}

\newcommand{\im}{{\Im}}
\newcommand{\re}{{\Re}}
\numberwithin{equation}{section}
\newtheorem{theorem}{Theorem}[section]
\newtheorem{proposition}[theorem]{Proposition}

\newtheorem{lemma}[theorem]{Lemma}
\newtheorem{definition}[theorem]{Definition}
\theoremstyle{definition}
\newtheorem{remark}[theorem]{Remark}

\newcommand{\brm}{\begin{remark}\rm}
\newcommand{\erm}{\end{remark}}
\newcommand{\brms}{\begin{remark}\rm}
\newcommand{\erms}{\end{remark}}
\newcommand{\bte}{\begin{theorem}}
\newcommand{\ete}{\end{theorem}}
\newcommand{\bpr}{\begin{proposition}}
\newcommand{\epr}{\end{proposition}}
\newcommand{\ble}{\begin{lemma}}
\newcommand{\ele}{\end{lemma}}
\newcommand{\beq}{\begin{equation}}
\newcommand{\eeq}{\end{equation}}
\newcommand{\bdm}{\begin{displaymath}}
\newcommand{\edm}{\end{displaymath}}
\numberwithin{equation}{section}

\newcommand{\bos}{\begin{remark}\rm}
\newcommand{\eos}{\end{remark}}
\newcommand{\dys}{\displaystyle}
\newcommand{\lbl}{\label}

\newcommand{\ben}{\begin{enumerate}}
\newcommand{\een}{\end{enumerate}}

\title[Soliton dynamics for NLS equations with magnetic field]
{Soliton dynamics for the nonlinear \\ Schr\"odinger equation with magnetic field}

\author[M.\ Squassina]{Marco Squassina }

\thanks{Department of Computer Science, 
University of Verona, 
C\`a Vignal 2, Strada Le Grazie 15, I-37134 Verona, Italy. 
E-mail: {\em marco.squassina@univr.it}}

\address{Department of Computer Science
\newline\indent
University of Verona
\newline\indent
C\'a Vignal 2, Strada Le Grazie 15, I-37134 Verona, Italy}
\email{marco.squassina@univr.it}

\thanks{The author was supported by the 2007 MIUR national research
project entitled {\em ``Variational and Topological Methods in the Study of
Nonlinear Phenomena''}}
\begin{document}

\begin{abstract}
The semiclassical regime of a nonlinear focusing Schr\"odinger equation in presence of non-constant 
electric and magnetic potentials $V,A$ is studied by taking as initial datum the ground state 
solution of an associated autonomous stationary equation. The concentration curve of the solutions is a parameterization
of the solutions of the second order ordinary equation $\ddot x=-\nabla V(x)-\dot x\times B(x)$, where $B=\nabla\times A$ is 
the magnetic field of a given magnetic potential $A$.
\end{abstract}
\maketitle


\bigskip
\begin{center}
\begin{minipage}{11cm}
	\footnotesize
\tableofcontents
\end{minipage}
\end{center}


\section{Introduction}
The aim of this paper is the study of the asymptotic behaviour of the solutions of the semilinear Schr\"odinger 
equation with an external magnetic potential $A$,
\begin{equation}
\label{probMF}
\tag{$P$}
\begin{cases}
\iu\eps\partial_t\phi_\eps=\frac{1}{2}\big(\frac{\eps}{\iu}\nabla-A(x)\big)^2\phi_\eps
+V(x)\phi_\eps-|\phi_\eps|^{2p}\phi_\eps,\quad & \text{$x\in\R^N,\, t>0$}, \\
\noalign{\vskip4pt}
\phi_\eps(x,0)=\phi_0(x),\quad & \text{$x\in\R^N$},
\end{cases}
\end{equation}
in the semiclassical regime of $\eps$ going to zero, by choosing a suitable class of initial data $\phi_0$ which is 
related to the (unique) ground state solution $r$ of an associated elliptic problem. We will show that
the evolution $\phi_\eps(t)$ remains close to $r$, in a suitable sense (and with an {\em explicit} convergence rate), locally uniformly in time,
in the transition from quantum to classical mechanics, namely as $\eps$ vanishes. 
This dynamical behaviour is also known as {\em soliton dynamics} (for a beautiful 
survey on {\em solitons} and their stability, see~\cite{tao-solitons}). Here, $\iu$ is the imaginary unit, $\eps$
is a small positive parameter playing the r\^ole of Planck's constant, $N\geq 1$, $0<p<2/N$ and
$V:\R^N\to\R$, $A:\R^N\to\R^N$ are an {\em electric} and {\em magnetic} potentials respectively. 
The magnetic field $B$ is $B=\nabla\times A$ in $\R^3$ and can be thought (and identified) in general dimension
as a $2$-form ${\mathbb H}^B$ of coefficients $(\partial_i A_j-\partial_j A_i)$. The magnetic Schr\"odinger operator which 
appears in problem~\eqref{probMF} formally operates as follows
\begin{equation}
\label{MagnOper}
\big(\frac{\eps}{\iu}\nabla-A(x)\big)^2\phi=-\eps^2\Delta\phi-\frac{2\eps}{\iu}A(x)\cdot\nabla\phi+|A(x)|^2\phi-\frac{\eps}{\iu}{\rm div}_xA(x)\phi,
\end{equation}
and it has been intensively studied in works by J.\ Avron, I.\ Herbst and B.\ Simon 
around 1978 (see~\cite{avron,avron2,avron3,ReedSimon,simonb} and references therein).
If $A=0$, then equation~\eqref{probMF} reduces to
\begin{equation}
	\label{ScA=0}
\begin{cases}
\iu\eps\partial_t\phi_\eps=-\frac{\eps^2}{2}\Delta\phi_\eps
+V(x)\phi_\eps-|\phi_\eps|^{2p}\phi_\eps,\quad & \text{$x\in\R^N,\, t>0$}, \\
\noalign{\vskip4pt}
\phi_\eps(x,0)=\phi_0(x),\quad & \text{$x\in\R^N$}.
\end{cases}
\end{equation}
For equation~\eqref{ScA=0}, rigorous results about the soliton dynamics were obtained in various papers
by J.C.\ Bronski, R.L.\ Jerrard \cite{bronski} and S.\ Keraani \cite{Keerani1,Keerani2} via 
arguments purely based on the use of conservation laws satisfied by the equation and by the associated
Hamiltonian system in $\R^N$ built upon the potential $V$, that is the {\em Newton law}
\begin{equation}
	\label{newtonsenza}
\ddot x=-\nabla V(x),\quad \dot x(0)=\xi_0,\,\, x(0)=x_0.
\end{equation}
For other achievements about the full dynamics of~\eqref{ScA=0}, see also~\cite{gril1,gril2} 
(in the framework of orbital stability of standing waves) as well as~\cite{kaup,keener} (in the framework
of non-integrable perturbation of integrable systems). Similar results were investigated in geometric optics 
by a different technique (WKB method), namely
writing formally the solution as $u_\eps=U_\eps(x,t) e^{\iu\theta(x,t)/\eps}$, where $U_\eps=U_0
+\eps U_1+\eps^2 U_2\cdots$, where $\theta$ and $U_j$ are solutions, respectively, of a Hamilton-Jacobi type
equation (known as {\em eikonal} equation) and of a system of transport equations.

It is very important to stress that, in the particular case of {\em standing wave solutions}  of~\eqref{ScA=0}, namely 
special solutions of~\eqref{ScA=0} of the form 
$$
\phi_\eps(x,t)=u_\eps(x) e^{-\frac{\iu}{\eps}\theta t},\quad x\in\R^N,\, t\in\R^+,\quad(\theta\in\R),
$$ 
where $u_\eps:\R^N\to\R$, there is an enormous
literature regarding the semiclassical limit for
the corresponding elliptic equation
$$
-\textstyle{\frac{\eps^2}{2}}\Delta u_\eps+V(x)u_\eps=|u_\eps|^{2p}u_\eps,\quad x\in\R^N.
$$
See the recent book~\cite{ambook} by A.\ Ambrosetti and A.\ Malchiodi and the references therein. Moreover,
there are various works admitting the presence of a magnetic potential $A$, and studying the asymptotic profile of the 
solutions $u_\eps:\R^N\to\C$ to 
$$
\textstyle{\frac{1}{2}}\big(\frac{\eps}{\iu}\nabla-A(x)\big)^2u_\eps+V(x)u_\eps=|u_\eps|^{2p}u_\eps,\quad x\in\R^N,
$$
as $\eps$ goes to zero (see e.g.~\cite{AS,BDP,CS,C,CS1,CSJJ,K,secsqu} and references therein).

In the special case $A=V=0$, the orbital stability for problem~\eqref{ScA=0} was proved by T.\ Cazenave and P.L.\ Lions~\cite{cl}, and by
M.\ Weinstein in~\cite{weinsteinMS,weinstein2}. Then, A.\ Soffer and M.\ Weinstein proved in~\cite{soffer1}
the asymptotic stability of nonlinear ground states of~\eqref{ScA=0}. 

See also the following seminal contributions ({\em in alphabetical order}): 
W.K.\ Abou Salem \cite{abou2},
V.\ Buslaev and G.\ Perelman~\cite{buslaev1,buslaev2},  
V.\ Buslaev and C.\ Sulem~\cite{buslaev3},
J.\ Fr\"ohlich, S.\ Gustafson, L.\ Jonsson, I.M.\ Sigal, T.-P.\ Tsai and H.-T.\ Yau~\cite{frolich1,frohl3,frohl4,frohl5},
S.\ Gustafson and M.I.\ Sigal~\cite{gustaf},
J.\ Holmer and Zworski~\cite{holmer,holmer2},
A.\ Soffer and M.\ Weinstein~\cite{soffer2,soffer3},
T.-P.\ Tsai and H.-T.\ Yau~\cite{tsai1,tsai2,tsai3}.
See also the references included in these works.

Now, in presence of a magnetic, some natural questions arise: what is the r\^ole played by the magnetic field $B$?
if $B$ plays a significant r\^ole, what is the correct Newton equation taking 
the place of~\eqref{newtonsenza}, which characterizes the 
concentrating curve and drives the dynamic in the semiclassical limit? 

As known, a charged particle moving in a magnetic field $B$ feels a sideways force that is proportional 
to the strength of $B$ as well as to its velocity. This force, which is always perpendicular to 
both the velocity of the particle and the magnetic field that created it (a particle moving in the direction of $B$ 
does not experience a force) is known as the {\em Lorentz force}. Hence, charged particles move in a circle (or more 
generally, {\em helix}) around the field lines of $B$ ({\em cyclotron motion}). During the motion, $B$ 
can do no work on a charged particle (cannot speed it up or slow it down) although it changes its direction 
(See figures~\ref{1fig} and~\ref{2fig}). 

As a consequence, with the expectation (which arises from the magnetic-free case) that in the semiclassical
limit the dynamics is governed by the classical Newtonian law, one is tempted to say that, in presence of an external magnetic field $B$, the right counterpart of~\eqref{newtonsenza} is given by the following Newton equation
\begin{equation}
	\label{NewtNew}
\ddot x=-\nabla V(x)-\dot x\times B(x),\quad \dot x(0)=\xi_0,\,\, x(0)=x_0,
\end{equation}
agreeing that $\times$ has to be interpreted as a matrix operation (${\mathbb H}^B\dot x$) if we are not in $\R^3$.

Only after full completion of the present paper the author discovered that a first result (mass and momentum asymptotics) in this direction
was obtained, independently, with decay assumptions on $B$ ,
by A.\ Selvitella in~\cite{selvit}, showing that, in fact, the above guess is the correct interpretation, in the transition process
from quantum to classical mechanics.

On the other hand, in this paper, we improve the result of~\cite{selvit}, proving a stronger result, which is precisely the one predicted by the WKB method.
Roughly speaking, under suitable regularity assumptions on $V$ and $A$, we show that, given the initial
position and velocity $x_0,\xi_0$ in $\R^N$, and taking as initial datum for~\eqref{probMF} 
\begin{equation}
\label{initialD}
\tag{$I$}
\phi_0(x)=r\Big(\frac{x-x_0}{\eps}\Big)
e^{\frac{\iu}{\eps}[A(x_0)\cdot(x-x_0)+x\cdot\xi_0]},\quad x\in\R^N,
\end{equation}
where $r\in H^1(\R^N)$ is the unique (up to translation) real 
ground state solution (bump like) of the associated elliptic problem
\begin{equation}
\label{seMF}
\tag{$S$}
-\frac{1}{2}\Delta r+r=|r|^{2p}r\qquad\text{in $\R^N$},
\end{equation}
then there exists a family of shift functions $\teta_{\vep}:\R^+\to [0,2\pi)$ such that
\begin{equation}
\label{mainconclusintro}
\phi_{\vep}(x,t)= r\Big(\frac{x-x(t)}{\vep}\Big)e^{\frac{{\rm i}}{\vep}\left[A(x(t))\cdot  (x-x(t))+x\cdot  \dot x(t)+\teta_{\vep}(t)\right]}+\omega_\eps,
\quad x\in\R^N,\,\,t>0,
\end{equation}
locally uniformly in time, as $\eps$ goes to zero, where we have set $\|\omega_\eps\|_{\H_\eps}={\mathcal O}(\eps)$,
and being $\|\phi\|_{\H_\eps}^2=\eps^{2-N}\|\nabla \phi\|_{L^2}^2+\eps^{-N}\|\phi\|_{L^2}^2$. In particular, with respect to~\cite{selvit},
the convergence rate is {\em explicit and of the order $\eps$} and, as a direct consequence,
the concentration center in the representation formula~\eqref{mainconclusintro} (expressing the soliton dynamics as guessed by 
the WKB method) is {\em exactly} $x(t)$ (in~\cite{selvit} formula~\eqref{mainconclusintro}
is not achievable, being the convergence rate undetermined).

The magnetic potential $A$ contributes to the {\em phase}
of the solution, and $x(t)$ is the {\em concentration line} (which can be considerably influenced by the presence of $B$,
see the phase portraits in figures~\ref{1fig}-\ref{2fig}).
Initial data~\eqref{initialD} should also be thought as corresponding to a {\em point particle}
with position $x_0$ and velocity $\xi_0$.
\vskip3pt

In the case where $\xi_0=0$ and $x_0$ is a {\em critical point} of the 
potential $V$, as equation~\eqref{NewtNew} admits the trivial solution $x(t)=x_0$ and $\xi(t)=0$
for all $t\in\R^+$, formula~\eqref{mainconclusintro} reduces to
\begin{equation*}
\phi_{\vep}(x,t)= r\Big(\frac{x-x_0}{\vep}\Big)
e^{\frac{{\rm i}}{\vep}\left[A(x_0)\cdot  (x-x_0)+\teta_{\vep}(t)\right]}+\omega_\eps,\quad x\in\R^N,\,\,t>0,
\end{equation*}
locally uniformly in time, as $\eps$ goes to zero (see Remark~\ref{criticpointV}). In turn, the concentration of $\phi_\eps$ is {\em static}
and takes place around the critical points of $V$, instead occurring along a smooth curve in $\R^N$. This is consistent
with the literature for the standing wave solutions mentioned above.

\vskip35pt
\begin{center}\textbf{Acknowledgements.}\end{center}
The author wish to thank Proff.\ I.M.\ Sigal, E.H.\ Lieb, T.-P.\ Tsai and S.\ Keraani for pointing out
useful information and bibliographic references related to the soliton dynamic of NLS equations. The author also thanks Proff.\ M.\ Degiovanni
and D.\ Fortunato for feedbacks on a preliminary version of the paper. The author was supported by the 2007 MIUR national research
project entitled: {\em ``Variational and Topological Methods in the Study of
Nonlinear Phenomena''}.

\vskip18pt
\begin{center}\textbf{Organization of the paper.}\end{center}
In Section~\ref{statementsect}, we introduce the functional framework, the tools and the 
ingredients needed to write the statement of the main result of the paper, Theorem~\ref{mainthBest}.\
In Section~\ref{prelminary-result}, we collect various preparatory results concerning
the characterization of the energy levels of the problem, in the semiclassical regime.\
In Section~\ref{approx-r}, we state the main approximation estimates for the solutions.\
In Section~\ref{massmom-id}, we get two integral identities for the evolution of the 
mass and momentum densities.\
In Section~\ref{massmomest}, we obtain the approximation results for the mass and momentum densities.
In Section~\ref{proofmain}, we obtain an error estimate. In turn,
we conclude the proof of the main result of the paper, Theorem~\ref{mainthBest}.\
Finally, In Section~\ref{conclusions}, we summarize the results obtained.

\vskip30pt
\begin{center}\textbf{Main notations.}\end{center}
\begin{enumerate}
\item The imaginary unit is denoted by $\iu$.
\item The complex conjugate of any number $z\in\C$ is denoted by $\bar z$.
\item The real part of a number $z\in\C$ is denoted by $\Re z$.
\item The imaginary part of a number $z\in\C$ is denoted by $\Im z$.
\item For all $z,w\in\C$ it holds $\Re(\bar z w)=\Re(z\bar w)$.
\item For all $z,w\in\C$ it holds $\Im(\bar z w)=-\Im(z\bar w)$.
\item The symbol $\R^+$ means the positive real line $[0,\infty)$.
\item The ordinary inner product between two vectors $a,b\in\R^N$ is denoted by $a \cdot b$.
\item The standard $L^2$ norm of a function $u$ is denoted by $\|u\|_{L^2}$.
\item The standard $L^\infty$ norm of a function $u$ is denoted by $\|u\|_{L^\infty}$.
\item The symbols $\partial_t$ and $\partial_j$ mean $\frac{\partial}{\partial t}$ and $\frac{\partial}{\partial x_j}$ respectively.
\item The symbol $\Delta$ means $\frac{\partial^2}{\partial x_{1}^2}+\cdots+\frac{\partial^2}{\partial x_{N}^2}$.
\item The symbol $C^m(\R^N)$, for $m\in\N$, denotes the space of functions with continuous derivatives up to
the order $m$. Sometimes $C^m(\R^N)$ is endowed with the norm
$$
\|\phi\|_{C^m}=\sum_{|\alpha|\leq m}\|D^\alpha\phi\|_{L^\infty}<\infty.
$$
\item The symbol $\int f$ stands for the integral of $f$ over $\R^N$ with the Lebesgue measure.
\item The symbol $C^{2*}$ denotes the dual space of $C^2$. The norm of a $\nu$ in $C^{2*}$ is
$$
\|\nu\|_{C^{2*}}=\sup\Big\{\Big|\int \phi(x)\nu dx\Big|:\,\phi\in C^2(\R^N),\,\,\|\phi\|_{C^2}\leq 1\Big\}.
$$
Clearly, $C^{2*}$ contains the space of bounded Radon measures.
\item $C$ denotes a generic positive constant, which may vary inside a chain of inequalities.
\item The first and second ordinary derivatives of $t\mapsto x(t)$ are denoted by $\dot x$ and $\ddot x$.
\item We use the Landau symbols. In particular ${\mathcal O}(\eps)$ is a generic function such that
the $\limsup$ of $\eps^{-1}{\mathcal O}(\eps)$ is finite, as $\eps$ goes to zero.
\end{enumerate}
\medskip

\section{Statement of the main result}
\label{statementsect}
\subsection{Functional setup and tools}

It is quite natural to consider operator~\eqref{MagnOper} on the Hilbert space $H_{A,\eps}$
defined by the closure of $C^\infty_c(\R^N;\C)$ under the scalar product
$$
(u,v)_{H_{A,\eps}}=\Re\int(D^\eps u\cdot \overline{D^\eps v}+V(x)u\overline{v})dx,
$$
where $D^\eps u=(D^\eps_1u,\dots,D^\eps_Nu)$ and $D_j^\eps=\iu^{-1}\eps \partial_j-A_j(x)$, with induced norm
$$
\|u\|_{H_{A,\eps}}^2=\int\Big|\frac{\eps}{\iu}\nabla u-A(x)u\Big|^2dx+\int V(x)|u|^2dx<\infty.
$$
The dual space of $H_{A,\eps}$ is denoted by $H_{A,\eps}'$, while the space $H_{A,\eps}^2$ is the set of $u$
such that
$$
\|u\|_{H_{A,\eps}^2}^2=\|u\|_{L^2}^2+\|\big(\frac{\eps}{\iu}\nabla-A(x)\big)^2u\|_{L^2}^2<\infty.
$$ 
Moreover, to problem~\eqref{probMF} it can be naturally associated the functional
$E:H_{A,\eps}\to\R$ (see also formula~\eqref{totenergy}) 
\begin{equation*}
E(u)=\frac{1}{2}\int\Big|\frac{\eps}{\iu}\nabla u-A(x)u\Big|^2dx+\int V(x)|u|^2dx-\frac{1}{p+1}\int |u|^{2p+2}dx.
\end{equation*}
Finally, we consider the functional 
${\mathcal E}:H^1(\R^N;\C)\to\R$ associated with~\eqref{seMF}
\begin{equation*}
{\mathcal E}(u)=\frac{1}{2}\int |\nabla u|^2dx-\frac{1}{p+1}\int |u|^{2p+2}dx.
\end{equation*}
It is a standard fact that the solution $r$ of~\eqref{seMF} is the unique (up to translation) solution
of the following minimization problem
\begin{equation}
\label{variatcaract-r}
{\mathcal E}(r)=\min\{{\mathcal E}(u): u\in H^1(\R^N),\,\|u\|_{L^2}=\|r\|_{L^2}\}.
\end{equation}
We also set
\begin{equation}
	\label{valueofm}
m:=\|r\|_{L^2}^2.
\end{equation}
Also, $r$ is radially symmetric and decreasing, belongs to $C^2(\R^N)\cap H^2(\R^N)$, and it decays exponentially 
together with its derivatives up to the order two, that is
\begin{equation}
\label{expdecay}
|D^\alpha r(x)|\leq Ce^{-\sigma |x|},\quad x\in\R^N,
\end{equation}
for some $\sigma,C>0$ and all $0\leq |\alpha|\leq 2$ (see e.g.~\cite{berlions1}).

\subsection{Well-posedness and conservation laws}

We recall that in~\cite[Section 9.1]{cazenave} (see also \cite{EL}), in the particular case
$N=3$ and when the external magnetic field $B=(b_1,b_2,b_3)$ is constant (thus $A$ is linear with respect to $x$), the (global) well-posedness
of problem~\eqref{probMF} in the (natural) energy space $H_{A,\eps}$ as well as the $H^2_{A,\eps}$-regularity 
of the flux for $H^2_{A,\eps}$-initial data was investigated (see Proposition~\ref{wellP} below) by T.\ Cazenave, M.\ Esteban and
P.L.\ Lions. Furthermore, 
in general dimension $N$ and for a general (smooth)
vector potential $A$, the (local) well-posedness in the energy space $H_{A,\eps}$  has been recently studied in~\cite{michel}
by L.\ Michel. We also wish to cite earlier papers by Y.\ Nakamura and A.\ Shimomura~\cite{naka1},
Y.\ Nakamura~\cite{naka2} as well as the important paper by K.\ Yajima~\cite{yaji}. In particular, in~\cite{naka1},
if $B$ has decay assumptions at infinity,
the problem is locally solved in the weighted space $\Sigma(2)\subset H^2(\R^N;\C)$ of functions $f$ in $L^2(\R^N;\C)$ such that
$\|x^\alpha D^\beta f\|_{L^2}<\infty$ for all $\alpha,\beta$ with $|\alpha|,|\beta|\geq 0$ and $0\leq |\alpha+\beta|\leq 2$
(notice that, via the decay~\eqref{expdecay}, the initial datum $\phi_0$ in~\eqref{initialD} belongs to the space $\Sigma(2)$).
Finally, see also \cite[Theorems 4.6.5 and 5.5.1]{cazenave} and an abstract result,
Lemma A.1, in the Appendix of~\cite{CazWei}, by T.\ Cazenave and F.B.\ Weissler. 

In order to prove the main result of the paper, we will assume (among other things)  that 
$A$ is globally bounded (together with its higher order derivatives). Clearly with this assumption the
well-posedness and regularity features for~\eqref{probMF} get easier to study. On the contrary, if $A$ is unbounded,
there are genuine regularity problems and the situation gets more involved \cite{pieroD}.

\begin{definition}
\label{admissibilty}
We say that a (sufficiently smooth) vector potential $A:\R^N\to\R^N$ is admissible with respect to problem~\eqref{probMF} if 
the following Proposition~\ref{wellP} holds for $A$.
\end{definition}

\begin{proposition}{\bf [well-posedness statement]}
\label{wellP}
Assume that $0<p<2/N$.\ Then, for every $\eps>0$ and all
$\phi_0\in\ H_{A,\eps}$, there exists a unique global solution
$$
\phi_\eps\in C(\R^+,H_{A,\eps})\cap C^1(\R^+,H_{A,\eps}')
$$ 
of problem~\eqref{probMF} with $\sup\limits_{t\in\R^+}\|\phi_\eps(t)\|_{H_{A,\eps}}<\infty$. 
Moreover, the mass associated with $\phi_\eps(t)$, 
\begin{equation*}
{\mathcal N}_{\vep}(t)=\frac1{\vep^{N}}\int |\phi_{\vep}(t)|^{2}dx,
\end{equation*}
as well as the total energy $E_\eps(t)=\eps^{-N}E(\phi_\eps(t))$ associated with~\eqref{probMF}
\begin{align}
	\label{totenergy}
 E_{\vep}(t)& =\frac{1}{2\vep^{N}}\int \left|\frac{\eps}{\iu}\nabla\phi_\eps(t)-A(x)\phi_{\vep}(t)\right|^{2}dx \\
&+\frac1{\vep^{N}}\int  V(x)|\phi_{\vep}(t)|^{2}dx 
-\frac{1}{(p+1)\eps^N}\int |\phi_{\vep}(t)|^{2p+2}dx, \notag
\end{align}
are conserved in time, namely 
$$
{\mathcal N}_{\vep}(t)={\mathcal N}_{\vep}(0),\quad
E_{\vep}(t)= E_{\vep}(0),\qquad\text{ for all $t\in\R^+$}.
$$
Finally if $\phi_0\in H_{A,\eps}^2$, then 
$\phi_\eps\in C(\R^+,H_{A,\eps}^2)\cap C^1(\R^+,L^2(\R^N;\C))$.
\end{proposition}

\begin{remark}
\label{remmassa}
From Proposition~\ref{wellP}, due to the choice of the initial data~\eqref{initialD}, the mass ${\mathcal N}_{\vep}(t)$ 
is also {\em independent} of $\vep$. Indeed, via the mass conservation and formula~\eqref{valueofm},
\beq\label{eqmi}
{\mathcal N_{\vep}(t)}={\mathcal N_{\vep}(0)}=\frac1{\vep^{N}}\int|\phi_{\vep}(x,0)|^{2}dx=\frac1{\vep^{N}}\int \Big|r\Big(\frac{x-x_0}{\vep}\Big)
\Big|^{2}dx=\|r\|^2_{L^2}=m,
\eeq
for all $\eps>0$ and $t\in\R^+$.
\end{remark}

\subsection{The driving Newtonian equation}
Given the initial data $x_0,\xi_0\in\R^N$, we consider 
$$
x(t),\xi(t):\R^+\to\R^{N},
$$ 
being the (unique) global (under the regularity assumptions on $V$ and $A$ indicated below) 
solution of the first order differential system
\begin{equation}
	\label{DriveS}
	\begin{cases}
		\dot x(t)=\xi(t), & \\
		\noalign{\vskip3pt}
		\dot \xi(t)=-\nabla V(x(t))-\xi(t)\times B(x(t)), & \\
		\noalign{\vskip1pt}
		x(0)=x_0, &\\
		\noalign{\vskip1pt}
        \xi(0)=\xi_0, &
	\end{cases}
\end{equation}
namely the second order ODE~\eqref{NewtNew}.
Notice that, setting
\begin{equation}
\label{Hamilt}
{\mathcal H}(t)=\frac{1}{2}|\xi(t)|^2+V(x(t)),\quad t\in\R^+,
\end{equation}
${\mathcal H}$ is a first integral associated with~\eqref{DriveS}, namely
$$
{\mathcal H}(t)={\mathcal H}(0),\qquad\text{for all $t\in\R^+$}.
$$ 
In general dimension $N$,
this follows by the elementary observation that, as ${\mathbb H}^B(x)$ is a skew-symmetric
matrix, we have $\xi(t)\cdot {\mathbb H}^B(x(t))\xi(t)=0$ for all $t\in\R^+$.
\medskip

\subsection{The main result}

We consider the following assumptions on the external electric and magnetic potentials, $V$ and $A$. 
\vskip6pt
\noindent
{\bf (V)} $V\in C^3(\R^N)$ is positive and $\|V\|_{C^3}<\infty$.
\vskip3pt
\noindent
{\bf (A)} $A\in C^3(\R^N;\R^N)$ with $\|A\|_{C^3}<\infty$ and $A$ is admissible (cf.\ definition~\ref{admissibilty}). 
\vskip6pt

Consider $H^1(\R^N;\C)$ equipped with the scaled norm $\|\phi\|_{\H_\eps}$,
$$
\|\phi\|_{\H_\eps}^2=\eps^{2-N}\|\nabla \phi\|_{L^2}^2+\eps^{-N}\|\phi\|_{L^2}^2.
$$
\vskip3pt
\noindent
The main result of the paper is the following

\begin{theorem}
\label{mainthBest}
Let $r$ be the ground state solution of problem~\eqref{seMF} and let $\phi_{\vep}$ be the family  
of solutions to problem~\eqref{probMF} with initial data~\eqref{initialD}, for some $x_0,\xi_0\in\R^N$.
Let $(x(t),\xi(t))$ be the global solution to system~\eqref{DriveS}. 
Then there exist $\delta>0$ and a locally uniformly bounded family of maps $\teta_{\vep}:\R^+\to [0,2\pi)$ 
such that, if $\|A\|_{C^2}<\delta$, then
\begin{equation}
\label{mainconclus}
\phi_{\vep}(x,t)=r\Big(\frac{x-x(t)}{\vep}\Big)
e^{\frac{{\rm i}}{\vep}\left[A(x(t))\cdot  (x-x(t))+x\cdot  \xi(t)+\teta_{\vep}(t)\right]}+\omega_\eps,   
\end{equation}
locally uniformly in time, where $\omega_\eps\in \H_{\vep}$ and 
$\|\omega_\eps\|_{\H_\vep}={\mathcal O}(\vep)$, as $\eps\to 0$.
Furthermore, without restrictions on $\|A\|_{C^2}$, we have
\begin{equation}
	\label{secondconcl-A}
	|\phi_\vep(x,t)|=r\Big(\frac{x-x(t)}{\vep}\Big)+\hat\omega_\eps^j(x,t),
\end{equation}
locally uniformly in time, where $\hat\omega_\eps\in \H_{\vep}$ and $\|\hat\omega_\eps^j\|_{\H_\vep}
\leq {\mathcal O}(\eps)$, as $\eps\to 0$.
\end{theorem}
\vskip4pt

Some comments are now in order.

\begin{remark}
	\label{criticpointV}
If $x_0$ is a {\em critical point} of $V$ and $\xi_0=0$, then the solution of system~\eqref{DriveS} is
$(x(t),\xi(t))=(x_0,0)$ for all $t\in\R^+$. Then, the conclusion of the previous result reads as
\begin{equation*}
\phi_{\vep}(x,t)=r\Big(\frac{x-x_0}{\vep}\Big)
e^{\frac{{\iu}}{\vep}\left[A(x_0)\cdot  (x-x_0)
+\teta_{\vep}(t)\right]}+\omega_\eps,
\end{equation*}
locally uniformly in time, where $\omega_\eps\in \H_{\vep}$ and $\|\omega_\eps\|_{\H_\vep}={\mathcal O}(\vep)$ as $\eps\to 0$. In particular, 
this is consistent with the literature of the {\em standing wave} solutions of~\eqref{probMF} in presence of a magnetic potential $A$ 
(see e.g.~\cite{AS,BDP,CS,C,CS1,CSJJ,K} and references included).
\end{remark}

\begin{remark}
In the framework of Theorem~\ref{mainthBest}, by the exponential decay of $r$, it holds
$$
|\phi_{\vep}(x,t)|\leq C e^{-\sigma\frac{|x-x(t)|}{\eps}}+|\omega_\eps(x,t)|.
$$
For an arbitrarily small $\delta>0$, the solution $\phi_\eps$ of~\eqref{probMF} is expected to decay exponentially 
in the set ${\mathcal P}_\delta=\{x\in\R^N:|x-x(t)|\geq\delta>0,\,\text{for all $t\in\R^+$}\}$ faster 
and faster as $\eps\to 0$, namely $\phi_\eps$ rapidly vanishes far from the concentration curve $x(t)$.
\end{remark}

\begin{remark}
A typical situation in $\R^3$ is when the external magnetic field $B=(b_1,b_2,b_3)$ is constant.
Without loss of generality, up to a rotation, one can assume that $B=(0,0,b)$ for some $b\in\R$.
Hence, the corresponding vector potential is $A(x,y,z)=\frac{b}{2}(-y,x,0)$. In this case, for harmonic external
potentials $V$, namely 
$$
V(x_1,x_2,x_3)=\frac{1}{2}(\omega_1^2x_1^2+\omega_2^2x_2^2+\omega_3^2x_3^2),\qquad\omega_j\in\R,
$$ 
system~\eqref{DriveS} reduces to
\begin{equation}
	\label{ex1}
	\begin{cases}
		\dot x_1(t)=\xi_1(t), & \\
		\dot x_2(t)=\xi_2(t), & \\
		\dot x_3(t)=\xi_3(t), & \\
		\dot \xi_1(t)=-\omega_1^2x_1(t)-b\xi_2(t), & \\
		\noalign{\vskip3pt}
	    \dot \xi_2(t)=-\omega_2^2x_2(t)+b\xi_1(t), & \\
		\noalign{\vskip3pt}
		\dot \xi_3(t)=-\omega_3^2x_3(t). & 
	\end{cases}
\end{equation}
It is clear that, setting some fixed values of $\omega_j$ and choosing some initial data,
enlarging the value of the third component $b$ of the magnetic field $B$ (say, from $0$ to $60$), the original periodic
orbit turns into a more and more helicoidal pattern. See figures~\ref{1fig} and~\ref{2fig}.
\begin{figure}[h!!!]
\begin{center}
   \includegraphics[scale=.78]{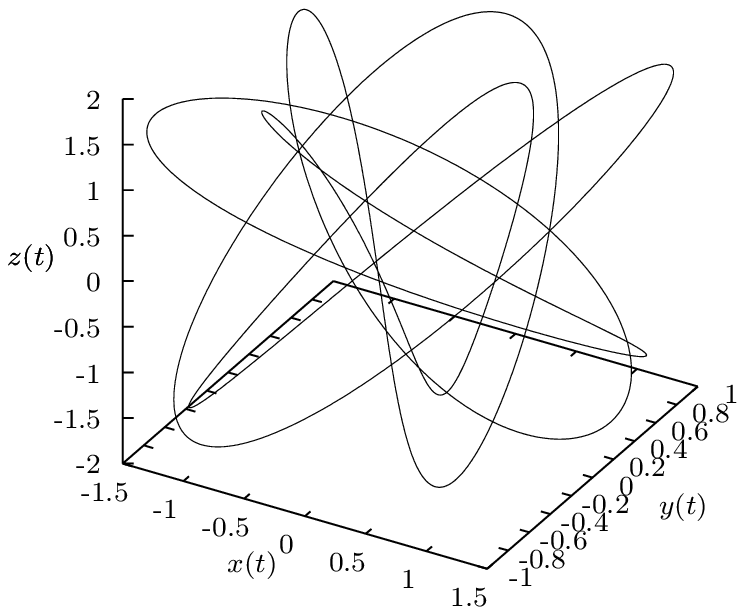}
   \hspace{20pt}
   \includegraphics[scale=.78]{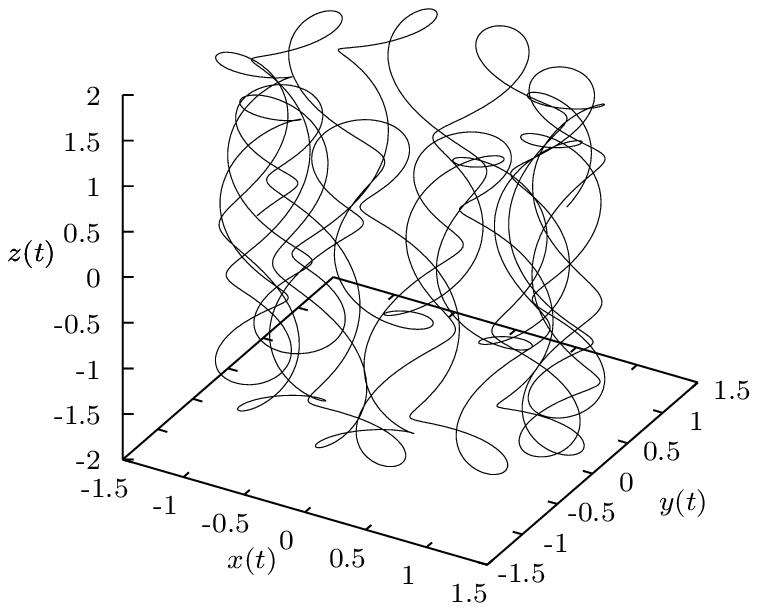}
\end{center}
\caption{Phase portrait of system~\eqref{ex1} with $\omega_1=1$, $\omega_2=1.4$, $\omega_3=1.2$, $b=0$~(left, no magnetic field) and 
	 $b=5$~(right, weak magnetic field).}
	 \label{1fig}
\end{figure}
\begin{figure}[h!!!]
\begin{center}
   \includegraphics[scale=.78]{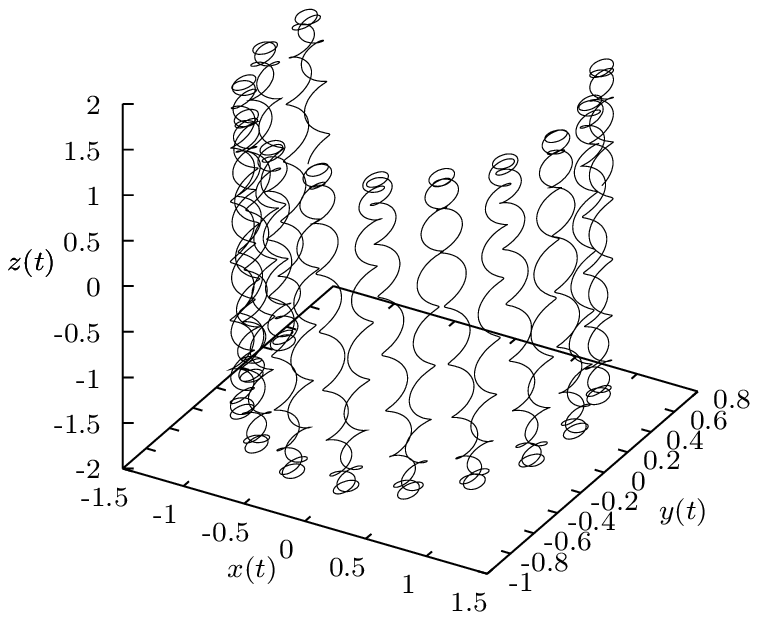}
   \hspace{20pt}
   \includegraphics[scale=.78]{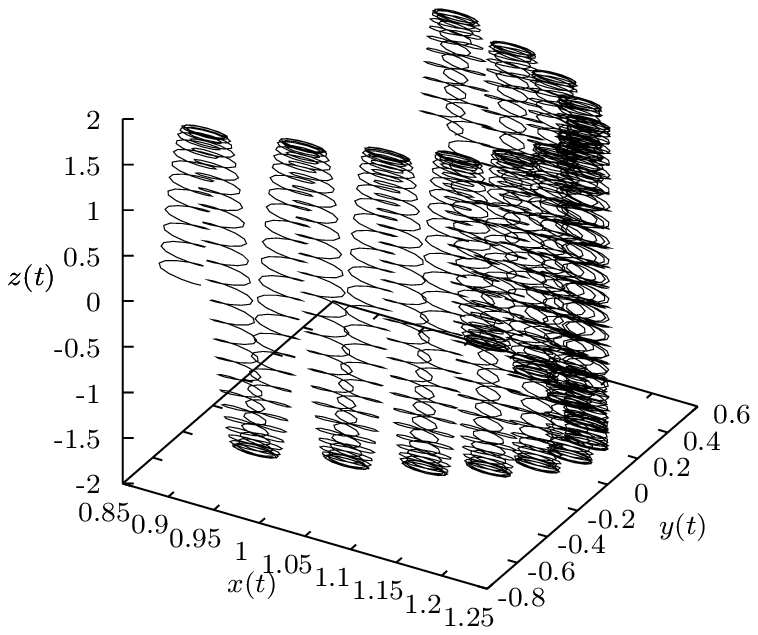}
\end{center}
\caption{Phase portrait of system~\eqref{ex1} with $\omega_1=1$, $\omega_2=1.4$, $\omega_3=1.2$, $b=20$~(left) and 
	 $b=60$~(right). The effects of the magnetic field increases.}
	 \label{2fig}
\end{figure}
\end{remark}	
\vskip2pt	

\begin{remark}
By complicating some arguments, assumption {\bf (V)} could be 
relaxed. For instance $V$ can be written as $V_1+V_2$, being $\|V_1\|_{C^3}<\infty$
and $V_2$ sufficiently smooth. The idea is to use the cut-off function indicated in~\eqref{defchi}, which 
is nonzero in the ball of $\R^N$ containing the region where the orbit $x(t)$ is confined (see~\cite{Keerani2}).
\end{remark}

\begin{remark}
As remarked in~\cite{Keerani2}, the soliton dynamics behaviour {\em breaks down in the critical case $p=\frac{2}{N}$}. Indeed, 
in this case, if we choose $x_0=\xi_0=0$, $V(x)=\frac{1}{2}|x|^2$ and $A=0$, then the modulus of the solution of problem~\eqref{probMF} 
with initial data $\phi_0(x)=r(\frac{x}{\eps})$ is given by $|\phi_\eps(x,t)|=(\cos t)^{-N/2}r(\frac{x}{\eps\cos t})$ for all $x\in\R^N$ and $t\in[0,\frac{\pi}{2})$
(see also \cite{c2}).
\end{remark}
\medskip

\section{Preliminary facts}
\label{prelminary-result}

In this section we collect some preliminary 
result which will allow us to prove the main result, Theorem~\ref{mainthBest}.

\subsection{Magnetic momentum}
The following vector function is useful to pursue our goals.

\begin{definition}
We define the momentum of the solution $\phi_\eps$, depending upon the vector potential $A$,
as a function $p^A_\eps:\R^N\times\R^+\to\R^N$, by setting 
\begin{equation*}
p^A_\eps(x,t):=\frac{1}{\eps^N}\im\big(\bar\phi_\eps(x,t)(\eps \nabla\phi_\eps(x,t)-\iu A(x)\phi_\eps(x,t)\big),\quad
x\in\R^N,\, t\in\R^+.
\end{equation*}
\end{definition}

First we state the following

\begin{lemma}
	\label{normBB}
Let $\phi_\eps$ be the solution to problem~\eqref{probMF} corresponding to the initial data~\eqref{initialD}.
Then there exists a positive constant $C$ such that
$$
\Big\|\frac{\eps}{\iu}\nabla\phi_\eps(\cdot,t)-A(x)\phi_{\vep}(\cdot,t)\Big\|^{2}_{L^2}\leq C\eps^N,
$$
for all $t\in\R^+$ and any $\eps>0$.
\end{lemma}
\begin{proof}
The total energy $E_\eps(t)$ is conserved (see Proposition~\ref{wellP}) and it
can be bounded independently of $\eps$ (see Lemma~\ref{estUNO}).
Then, since $V$ is positive, defining $\zeta_\eps(x):=\phi_\eps(\eps x)$, 
it follows that, for some positive constant $C$,
\begin{equation}
	\label{first}
\Big\|\frac{1}{\iu}\nabla\zeta_\eps(\cdot,t)-A(\eps x)\zeta_{\vep}(\cdot,t)\Big\|^{2}_{L^2}
-C\|\zeta_{\vep}(\cdot,t)\|^{2p+2}_{L^{2p+2}}\leq C.
\end{equation}
By combining the diamagnetic inequality (see \cite{EL} for a proof)
$$
|\nabla |\zeta_\eps||\leq \Big|\Big(\frac{\nabla}{\iu}-A(\eps x)\Big)\zeta_\eps\Big|,\qquad \text{a.e.\ in $\R^N$}
$$
with the Gagliardo-Nirenberg inequality, setting $\vartheta=\frac{pN}{2p+2}\in (0,1)$, we obtain 
$$
\|\zeta_{\vep}(\cdot,t)\|_{L^{2p+2}}\leq \|\zeta_{\vep}(\cdot,t)\|_{L^{2}}^{1-\vartheta}
\|\nabla |\zeta_{\vep}(\cdot,t)|\|^{\vartheta}_{L^{2}}\leq \|\zeta_{\vep}(\cdot,t)\|_{L^{2}}^{1-\vartheta}
\Big\|\Big(\frac{\nabla}{\iu}-A(\eps x)\Big)\zeta_\eps(\cdot,t)\Big\|^{\vartheta}_{L^{2}}.
$$
By the conservation of mass (see Remark~\ref{remmassa}), we deduce that $\|\zeta_\eps(\cdot,t)\|_{L^2}^2
={\mathcal N}_\eps(t)=m$, independently of $\eps$. Hence,  for all $\eps>0$, we get
$$
\|\zeta_{\vep}(\cdot,t)\|_{L^{2p+2}}^{2p+2}
\leq C\Big\|\frac{1}{\iu}\nabla\zeta_\eps(\cdot,t)-A(\eps x)\zeta_\eps(\cdot,t)\Big\|^{pN}_{L^{2}}.
$$
Since $pN<2$ by assumption, the assertion readily follows from~\eqref{first} and rescaling.
\end{proof}

We have the following summability property for $p^A_\eps(x,t)$.

\begin{lemma}
	\label{boundedMoment}
	There exists a positive constant $C$ such that
	$$
	\sup_{t\in\R^+}\left|\int p^A_\eps(x,t)dx\right|\leq C.
	$$
\end{lemma}
\begin{proof}
Taking into account the inequality of Lemma~\ref{normBB} and the mass conservation law, 
by H\"older inequality we get
\begin{align*}
\left|\int p^A_\eps(x,t)dx\right| &\leq \int |p^A_\eps(x,t)|dx\leq
\frac{1}{\eps^N}\int |\bar\phi_\eps(x,t)|\Big|\frac{\eps}{\iu} \nabla\phi_\eps(x,t)-A(x)\phi_\eps(x,t)\Big|dx	\\
&\leq \frac{1}{\eps^{N/2}}\|\phi_\eps(\cdot,t)\|_{L^2}\frac{1}{\eps^{N/2}}\Big\|\frac{\eps}{\iu} \nabla\phi_\eps(\cdot,t)-A(x)\phi_\eps(\cdot,t)\Big\|_{L^2}\leq C,
\end{align*}
for all $t\in\R^+$. The assertion follows by taking the supremum over positive times.
\end{proof}

\subsection{Energy levels in the semiclassical limit}

Let us recall a useful tool (see e.g.~\cite[Lemma 3.3]{Keerani2}), which reveals useful 
in managing various estimates that follow.

\begin{lemma}\label{pote}
Assume that $g:\R^N\to\R$ is a function of class $C^{2}(\R^N)$, $\|g\|_{C^2}<\infty$, and that $r$ 
is the ground state solution of~\eqref{seMF}. Then, as $\eps$ goes to zero, it holds
\begin{equation*}
\int g(\vep x+y)r^2(x)dx=\int g(y)r^{2}(x)dx+{\mathcal O}(\eps^2),
\end{equation*}
for every $y\in \R^N$ fixed. Moreover, ${\mathcal O}(\eps^2)$ is uniform with respect
to a family ${\mathcal F}\subset C^2(\R^N)$ which is uniformly bounded, that is
$\sup_{g\in{\mathcal F}}\|g\|_{C^2}<\infty$.
\end{lemma}

In the next lemma we compute the value of the energy associated 
with~\eqref{probMF}-\eqref{initialD}, in the semiclassical regime.

\begin{lemma}
	\label{estUNO}
Let $E_\eps$ be the energy associated with the family $\phi_{\vep}$ of solutions to problem~\eqref{probMF} 
with initial data~\eqref{initialD}. Then, for every $t\in\R^+$, it holds
$$
E_\eps(t)={\mathcal E}(r)+m{\mathcal H}(t)+{\mathcal O}(\eps^2),
$$
as $\eps$ goes to zero.
\end{lemma}
\begin{proof}
Notice that, for all $x\in\R^N$, we get
\begin{align*}
& \left(\frac{\eps}{\iu}\nabla-A(x)\right)\Big(r\Big(\frac{x-x_0}{\eps}\Big)e^{\frac{\iu}{\eps}[A(x_0)\cdot(x- x_0)+x\cdot\xi_0]}\Big)=
\frac{1}{\iu}e^{\frac{\iu}{\eps}[A(x_0)\cdot(x-x_0)+x\cdot\xi_0]}\nabla r\Big(\frac{x-x_0}{\eps}\Big) \\
& +r\Big(\frac{x-x_0}{\eps}\Big)e^{\frac{\iu}{\eps}[A(x_0)\cdot(x-x_0)+x\cdot\xi_0]}[A(x_0)+\xi_0]-
r\Big(\frac{x-x_0}{\eps}\Big)e^{\frac{\iu}{\eps}[A(x_0)\cdot(x-x_0)+x\cdot\xi_0]}A(x).
\end{align*}
Hence, it follows that
\begin{align*}
& \frac{1}{\vep^{N}}\int\Big|\Big(\frac{\eps}{\iu}\nabla-A(x)\Big)
\Big(r\Big(\frac{x-x_0}{\eps}\Big)e^{\frac{\iu}{\eps}[A(x_0)\cdot(x-x_0)+x\cdot\xi_0]}\Big)\Big|^2dx=
\frac{1}{\vep^{N}}\int\Big|\nabla r\Big(\frac{x-x_0}{\eps}\Big)\Big|^2dx \\
& +\frac{1}{\vep^{N}}\int r^2\Big(\frac{x-x_0}{\eps}\Big)|A(x_0)+\xi_0|^2dx+
\frac{1}{\vep^{N}}\int r^2\Big(\frac{x-x_0}{\eps}\Big)|A(x)|^2dx  \\
& -\frac{2}{\vep^{N}}\int r^2\Big(\frac{x-x_0}{\eps}\Big)A(x)\cdot (A(x_0)+\xi_0)dx \\
&=\int\left|\nabla r(x)\right|^2dx +|A(x_0)+\xi_0|^2m+
\int r^2(x)|A(\eps x+x_0)|^2dx  \\
& -2\int r^2(x)A(\eps x+x_0)\cdot (A(x_0)+\xi_0)dx.
\end{align*}
In view of Lemma~\ref{pote}, we have
\begin{align*}
& \int r^2(x)|A(\eps x+x_0)|^2dx=|A(x_0)|^2m+{\mathcal O}(\eps^2),   \\
& \int r^2(x)A(\eps x+x_0)\cdot (A(x_0)+\xi_0)dx=A(x_0)\cdot (A(x_0)+\xi_0)m+{\mathcal O}(\eps^2).
\end{align*}
Then,
\begin{align*}
& \frac{1}{\vep^{N}}\int\Big|\Big(\frac{\eps}{\iu}\nabla-A(x)\Big)
\Big(r\Big(\frac{x-x_0}{\eps}\Big)e^{\frac{\iu}{\eps}[A(x_0)\cdot(x-x_0)+x\cdot\xi_0]}\Big)\Big|^2dx \\
&=\int\left|\nabla r(x)\right|^2dx +|A(x_0)+\xi_0|^2m+|A(x_0)|^2m  
-2A(x_0)\cdot (A(x_0)+\xi_0)m+{\mathcal O}(\eps^2) \\
&=\int\left|\nabla r(x)\right|^2dx +m|\xi_0|^2+{\mathcal O}(\eps^2).
\end{align*}
It turn, by combining the conservation of energy (see Proposition~\ref{wellP}) and the conservation 
of the function ${\mathcal H}$ (see definition~\eqref{Hamilt}), we get
\begin{align*}
E_\eps(t)&=E_\eps(0) =E_\eps\Big(r\Big(\frac{x-x_0}{\eps}\Big)e^{\frac{\iu}{\eps}[A(x_0)\cdot(x-x_0)
+x\cdot\xi_0]}\Big) \\
&=\frac{1}{2\vep^{N}}\int\Big|\Big(\frac{\eps}{\iu}\nabla-A(x)\Big)
\Big(r\Big(\frac{x-x_0}{\eps}\Big)e^{\frac{\iu}{\eps}[A(x_0)\cdot(x-x_0)+x\cdot\xi_0]}\Big)\Big|^2dx \\
& +\int  V(x_0+\eps x)r^{2}(x)dx -\frac{1}{p+1}\int |r(x)|^{2p+2}dx \\
& =\frac{1}{2}\int\left|\nabla r(x)\right|^2dx +\int  V(x_0+\eps x)r^{2}(x)dx -\frac{1}{p+1}\int |r(x)|^{2p+2}dx
+\frac{1}{2}m|\xi_0|^2+{\mathcal O}(\eps^2) \\
&={\mathcal E}(r)+\int  V(x_0+\eps x)r^{2}(x)dx+\frac{1}{2}m|\xi_0|^2+{\mathcal O}(\eps^2) \\
& ={\mathcal E}(r)+mV(x_0)+\frac{1}{2}m|\xi_0|^2+{\mathcal O}(\eps^2) \\
& ={\mathcal E}(r)+m{\mathcal H}(0)+{\mathcal O}(\eps^2) \\
\noalign{\vskip3pt}
& ={\mathcal E}(r)+m{\mathcal H}(t)+{\mathcal O}(\eps^2),
\end{align*}
as $\eps$ goes to zero.
\end{proof}

\begin{lemma}
	\label{estDUE}
Let  $\phi_{\vep}$ be the family of solutions to problem~\eqref{probMF} with initial data~\eqref{initialD}.
Let us set, for any $\eps>0$, $t\in\R^+$ and $x\in\R^N$
\begin{equation}
\label{veps}
\psi_{\vep}(x,t) =e^{-\frac{{\rm i}}{\vep}\xi(t)\cdot [\vep x+x(t)]} e^{-\iu A(x(t))\cdot x}
\,\phi_{\vep}(\vep x+x(t),t)
\end{equation}
where $(x(t), \xi(t))$ is the solution of system~\eqref{DriveS}. Then
\begin{align*}
\E(\psi_{\vep}(t)) & = E_{\vep}(t) -\int V(x)\frac{|\phi_{\vep}(x,t)|^{2}}{\eps^N}dx
+\frac12 m|\xi(t)+A(x(t))|^2 \\
& -(\xi(t)+A(x(t))\cdot \int p_{\vep}^A(x,t)dx-(\xi(t)+A(x(t))\cdot \int A(x)\frac{|\phi_\eps(x,t)|^2}{\eps^N}dx\\
& +\frac{1}{2}\int |A(x)|^2\frac{|\phi_\eps(x,t)|^2}{\eps^N}dx
+\int A(x)\cdot p_\eps^A(x,t)dx.
\end{align*}
\end{lemma}
\begin{proof}
By a simple change of variable and Remark~\ref{remmassa}, we have  
\begin{equation}\label{normel2}
\|\psi_{\vep}(t)\|_{L^2}^{2}=\dys
\|\phi_{\vep}(\vep x+x(t),t)\|_{L^2}^{2}=
\frac1{\vep^{N}}\|\phi_{\vep}(t)\|_{L^2}^{2}={\mathcal N}_\eps(t)=m,\quad t\in\R^+.
\end{equation}
Hence the mass of $\psi_{\vep}(t)$ is conserved during the evolution. Let
\begin{equation*}
p_\eps(x,t)=\frac{1}{\eps^{N-1}}\im\big(\bar\phi_\eps(x,t)\nabla\phi_\eps(x,t)\big),\quad
x\in\R^N,\, t\in\R^+,
\end{equation*}
be the magnetic-free momentum. A direct computation yields
\begin{align*}
&\E(\psi_{\vep}(t))=\frac1{2\vep^{N-2}}\int |\nabla \phi_{\vep}(t)|^{2}dx +
\frac12 m|\xi(t)+A(x(t))|^2    \\
&-\frac1{\vep^{N}}\frac{1}{p+1}\int |\phi_{\vep}(t)|^{2p+2}dx 
-(\xi(t)+A(x(t))\cdot \int p_{\vep}(x,t)dx \\
& =\frac1{2\vep^{N}}\int\left|\frac{\eps}{\iu}\nabla \phi_{\vep}(t)-A(x)\phi_\eps(t)\right|^{2}dx \\
&-\frac{1}{2\eps^N}\int |A(x)|^2|\phi_\eps(t)|^2 dx+\frac{1}{\eps^{N-1}}\int A(x)\cdot 
\im(\bar\phi_\eps(t) \nabla\phi_\eps(t)) \\
& +\frac12 m|\xi(t)+A(x(t))|^2-\frac1{\vep^{N}}\frac{1}{p+2}\int |\phi_{\vep}(t)|^{2p+2}dx 
-(\xi(t)+A(x(t))\cdot \int p_{\vep}(x,t)dx. 
\end{align*}
Then, taking into account the definition of $E_\eps(t)$, we obtain
\begin{align*}
\E(\psi_{\vep}(t)) & = E_{\vep}(t) -\int V(x)\frac{|\phi_{\vep}(x,t)|^{2}}{\eps^N}dx+\frac12 m|\xi(t)+A(x(t))|^2   \\
&-(\xi(t)+A(x(t))\cdot \int p_{\vep}(x,t)dx  \\
& -\frac{1}{2\eps^N}\int |A(x)|^2|\phi_\eps(x,t)|^2dx+\int A(x)\cdot p_\eps(x,t)dx.
\end{align*}
Finally, since 
$$
p_\eps(x,t)=p_\eps^A(x,t)+\eps^{-N}A(x)|\phi_\eps(x,t)|^2,
$$ 
we obtain the desired conclusion.
\end{proof}
\vskip3pt

Next we introduce two important functionals in the dual space of $C^2$.

\begin{definition}
	\label{2dualdefn}
	Let  $\phi_{\vep}$ be the family of solutions to problem~\eqref{probMF} with initial data~\eqref{initialD}
	and let $p^A_\eps$ be the corresponding momentum.
For any $t\in\R^+$, let us define an element $\Pi^1_\eps(\cdot,t)$ in the dual space of $C^2(\R^N;\R^N)$ and
an element $\Pi^2_\eps(\cdot,t)$ in the dual space of $C^2(\R^N;\R)$ by setting
\begin{align*}
\forall\varphi\in C^2(\R^N;\R^N):\quad
&\int \Pi^1_\eps(x,t) \cdot\varphi\, dx=\int \varphi\cdot p^A_\eps(x,t)dx-m\varphi(x(t))\cdot\xi(t),  \\
\forall\varphi\in C^2(\R^N;\R):\quad
&\int \Pi^2_\eps(x,t)\varphi \,dx=\int \varphi\frac{|\phi_\eps(x,t)|^2}{\eps^N}dx-m\varphi(x(t)),
\end{align*}
and all $t\in\R^+$. Here $x(t),\xi(t)$ denote the components of the solution of system~\eqref{DriveS}.
\end{definition}

We recall a property of the functional $\delta_{y}$ on $C^{2}(\rn)$ (see~\cite[Lemma 3.1, 3.2]{Keerani2}).

\begin{lemma}\lbl{propdelta}
There exist three positive constants $K_{0},\,K_{1},\,K_{2}$ such that, for all $y,z\in\R^N$,
$$
K_{1}|y-z|\leq \|\delta_{y}-\delta_{z}\|_{C^{2*}}\leq K_{2}|y-z|,
$$
provided that $\|\delta_{y}-\delta_{z}\|_{C^{2*}}\leq K_{0}$.
\end{lemma}

For a fixed time $T_0>0$ (to be chosen later on), let $\rho$ be a positive constant defined by
\begin{equation}
	\label{defa}
\rho=K_{1}\sup_{[0,T_{0}]}|x(t)|+K_{0}
\end{equation}
where $x(t)$ is defined in~\eqref{DriveS}, the constants $K_{0}$ and $K_1$
are defined in Lemma~\ref{propdelta}, and  let $\chi$ be a 
$C^{\infty}(\rn)$ function such that $0\leq\chi\leq 1$ and 
\begin{equation}
	\lbl{defchi}
\chi(x)=1\quad\text{if $|x|< \rho$},\qquad\quad \chi(x)=0\quad
\text{if $|x|> 2\rho$}. 
\end{equation}

\noindent
Let us now set, for all $t\in\R^+$ and $\eps>0$,
\begin{align*}
\omega_\eps^1(t)&:=\int (\xi(t)+A(x(t))\cdot \Pi^1_\eps(x,t) dx, \\
\omega^2_\eps(t)&:=\int A(x)\cdot \Pi^1_\eps(x,t) dx,   \\
\omega^3_\eps(t)&:=\int |A(x)|^2\Pi_\eps^2(x,t) dx, \\
\omega^4_\eps(t)&:=\int (\xi(t)+A(x(t))\cdot A(x)\Pi_\eps^2(x,t)dx,  \\ 
\omega^5_\eps(t)&:=\int V(x)\Pi_\eps^2(x,t) dx, \\ 
\gamma_{\vep}(t)&:=mx(t)-\int x\chi(x)\frac{|\phi_{\vep}(x,t)|^{2}}{\eps^N}dx,
\end{align*}
where $\chi$ is as in~\eqref{defchi}.
\vskip5pt

On the functions $\omega^j_\eps$, we have the following estimate.
\begin{lemma}
	\label{controestimateee}
There exists a positive constant $C=C(V,A)$ such that 
\begin{equation}
	\label{dualcontrol}
	\sum_{j=1}^5|\omega^j_\eps(t)|\leq C\Omega_\eps(t),
\end{equation}
where the function $\Omega_\eps:\R^+\to\R^+$ is defined as 
$\Omega_\eps(t)=\hat\Omega_\eps(t)+\rho^A_\eps(t)$, where
\begin{align}
\label{defOmega}
\hat\Omega_\eps(t) &:=\left|\int \Pi^1_\eps(x,t)dx\right|+\sup_{\|\varphi\|_{C^3}\leq 1}
\left|\int \varphi\Pi^2_\eps(x,t)dx\right|+|\gamma_\eps(t)|,\quad t\in\R^+, \\
\label{defrho}
\rho^A_\eps(t)&:=\left|\int A(x)\cdot \Pi^1_\eps(x,t) dx\right|,\quad t\in\R^+.
\end{align}
Moreover 
$$
\Omega_\eps(0)={\mathcal O}(\eps^2),
$$ 
as $\eps$ goes to zero.
\end{lemma}
\begin{proof}
Estimate~\eqref{dualcontrol} is a simple and direct consequence of the definition of $\omega^j_\eps(t)$, $\Omega_\eps(t)$,
of the uniform boundedness of $\xi(t),A(x(t))$, namely $|\xi(t)|+|A(x(t))|\leq C$ and of the
fact that $\|V\|_{C^3}<\infty$ and $\|A\|_{C^3}<\infty$. Let us now prove that 
$\Omega_\eps(0)={\mathcal O}(\eps^2)$, as $\eps\to 0$. Recalling that the 
initial data $\phi_0$ is $r((x-x_0)/\eps)e^{\iu/\eps[A(x_0)\cdot(x-x_0)+x\cdot \xi_0]}$, 
in light of Lemma~\ref{pote}, for any $\varphi\in C^2(\R^N;\R^N)$ such that $\|\varphi\|_{C^2}\leq 1$,
we infer 
\begin{align*}
	& \int \varphi(x)\cdot \Pi^1_\eps(x,0) dx =\int \varphi(x)\cdot p^A_\eps(x,0)dx-m\varphi(x_0)\cdot \xi_0 \\  
	 & =\frac{1}{\eps^{N-1}}\int \varphi(x)\cdot\im\big(\bar\phi_\eps(x,0)\nabla\phi_\eps(x,0)\big)    \\
&	-\frac{1}{\eps^N}\int \varphi(x)\cdot A(x)|\phi_\eps(x,0)|^2dx-m\varphi(x_0)\cdot \xi_0 \\
	& =\frac{1}{\eps^{N}}\int \varphi(x)\cdot (A(x_0)+\xi_0) r^2\left(\frac{x-x_0}{\eps}\right)dx    \\          
	&-\frac{1}{\eps^{N}}\int \varphi(x)\cdot A(x)r^2\left(\frac{x-x_0}{\eps}\right)dx-m\varphi(x_0)\cdot \xi_0 \\
		& =\int \varphi(x_0+\eps x)\cdot(A(x_0)+\xi_0)r^2(x)dx   \\             
	& -\int \varphi(x_0+\eps x)\cdot A(x_0+\eps x)r^2(x)dx-m\varphi(x_0)\cdot \xi_0 \\
	& =m\varphi(x_0)\cdot(A(x_0)+\xi_0)                
	-m\varphi(x_0)\cdot A(x_0)    \\
	\noalign{\vskip4pt}
	&-m\varphi(x_0)\cdot \xi_0+{\mathcal O}(\eps^2)={\mathcal O}(\eps^2),
\end{align*}
as $\eps$ goes to zero. In a similar fashion, for any $\varphi\in C^3(\R^N)$ with $\|\varphi\|_{C^3}\leq 1$, we get
\begin{align*}
	\int \varphi(x)\Pi^2_\eps(x,0) dx & =\frac{1}{\eps^N}\int \varphi(x) |\phi_\eps(x,0)|^2 dx-m\varphi(x_0) \\
	& =\int \varphi(x_0+\eps x) r^2(x)dx-m\varphi(x_0)={\mathcal O}(\eps^2).
\end{align*}
Finally,  as $\chi(x_0)=1$, we have 
$|\gamma_{\vep}(0)|=\big|m x_0-\int (x_0+\vep y)\chi(x_0+\vep y)r^2(y)dy\big|\leq{\mathcal O}(\eps^2)$,
by  Lemma~\ref{pote}. This concludes the proof of the assertion.
\end{proof}

At this stage, we are ready to estimate the energy values $\E(\psi_{\vep}(t))$.

\begin{lemma}
	\label{firstKey}
Let $\psi_\eps$ be the function defined in formula~\eqref{veps}. Then there exists a positive constant $C$ such that
\begin{equation*}
0\leq \E(\psi_{\vep}(t)) -\E(r) \leq C\Omega_\eps(t)+{\mathcal O}(\eps^2),
\end{equation*}
for all $t\in\R^+$ and $\eps>0$.
\end{lemma}
\begin{proof}
By combining the conclusions of Lemma~\ref{estUNO} and~\ref{estDUE}, we obtain
\begin{align*}
\E(\psi_{\vep}(t)) -\E(r) & = m{\mathcal H}(t) -\int V(x)\frac{|\phi_{\vep}(x,t)|^{2}}{\eps^N}dx
+\frac12 m|\xi(t)+A(x(t))|^2 \\
& -(\xi(t)+A(x(t)))\cdot \int p_{\vep}^A(x,t)dx-(\xi(t)+A(x(t))\cdot \int A(x)\frac{|\phi_\eps(x,t)|^2}{\eps^N}dx\\
& +\frac{1}{2}\int |A(x)|^2\frac{|\phi_\eps(x,t)|^2}{\eps^N}dx
+\int A(x)\cdot p_\eps^A(x,t)dx+{\mathcal O}(\eps^2),
\end{align*}
for all $t\in \R^+$, as $\eps$ goes to zero. Notice that
\begin{align*}
	(\xi(t)+A(x(t))\cdot \int p_{\vep}^A(x,t)dx&=m|\xi(t)|^2+mA(x(t))\cdot\xi(t)+\omega_1(t),  \\
	\int A(x)\cdot p_\eps^A(x,t)dx&=mA(x(t))\cdot\xi(t)+\omega_2(t),   \\
	\int |A(x)|^2\frac{|\phi_\eps(x,t)|^2}{\eps^N}dx&=m|A(x(t))|^2+\omega_3(t), \\
 \int (\xi(t)+A(x(t))\cdot A(x)\frac{|\phi_\eps(x,t)|^2}{\eps^N}dx&=m\xi(t)\cdot A(x(t))+m|A(x(t))|^2+\omega_4(t), \\
\int V(x)\frac{|\phi_\eps(x,t)|^2}{\eps^N}dx&=m V(x(t))+\omega_5(t).
\end{align*}
It follows that
\begin{align*}
\E(\psi_{\vep}(t)) -\E(r) & = \frac{1}{2}m|\xi(t)|^2+mV(x(t))-mV(x(t))-\omega_5(t)
+\frac12 m|\xi(t)+A(x(t))|^2 \\
\noalign{\vskip3pt}
& -m|\xi(t)|^2-mA(x(t))\cdot\xi(t)-\omega_1(t)-m\xi(t)\cdot A(x(t))-m|A(x(t))|^2-\omega_4(t)  \\
\noalign{\vskip3pt}
& +\frac{1}{2}m|A(x(t))|^2+\frac{\omega_3(t)}{2}+mA(x(t))\cdot\xi(t)+\omega_2(t)\\
& =-\omega_1(t)+\omega_2(t) +\frac{\omega_3(t)}{2}-\omega_4(t)-\omega_5(t)+{\mathcal O}(\eps^2),
\end{align*}
which concludes the proof in light of inequality~\eqref{dualcontrol}
of Lemma~\ref{controestimateee}.
\end{proof}

\section{The approximation result}
\label{approx-r}

Let us first recall a useful and well-established stability property of ground states.

\begin{proposition}\label{stabilita}
There exist two positive constants $\cala$ and ${\mathcal C}$ such that, 
if $\Phi\in H^1(\R^N;\C)$ is such that $\|\Phi\|_{L^2}=\|r\|_{L^2}$, 
where $r$ is the ground state solution of~\eqref{seMF}, and
$$
\E(\Phi)-\E(r)\leq \cala,
$$
then 
\beq
\inf_{y\in\R^N,\,\vartheta\in [0,2\pi)}\|\Phi-e^{\iu \theta}r(\cdot +y)\|_{H^1}^{2}\leq {\mathcal C}\left(\E(\Phi)-\E(r)\right).
\eeq
\end{proposition}

\begin{proof}
See~\cite{weinsteinMS,weinstein2}.
\end{proof}

Next, in view of the previous preparatory work, we can state the reppresentation result.

\begin{theorem}
\label{chiave}
Let  $\phi_{\vep}$ be the family of solutions to problem~\eqref{probMF} with initial data~\eqref{initialD}
and let $\psi_{\eps}$ be the function defined in formula~\eqref{veps}.
Then there exist $\eps_{0}>0$, a time $T_{\eps}^{*}>0$, families of uniformly 
bounded functions $\theta_{\vep}:\R^+\to[0,2\pi)$, 
$y_{\vep}:\R^+\to\R^N$ and a positive constant $C$ such that
\begin{equation}
	\lbl{quasifine}
	\phi_\eps(x,t)=e^{\frac{\iu}{\vep}(\xi(t)\cdot x+\theta_\eps(t)+A(x(t))
	\cdot (x-x(t))}r\Big(\frac{x-y_\vep(t)}{\vep}\Big)+\omega_\eps(t),
\end{equation}	
where
\begin{equation*}
	\|\omega_\eps(t)\|_{\H_{\vep}}\leq C\sqrt{\Omega_\eps(t)}+{\mathcal O}(\eps),
\end{equation*}
for all $\eps\in(0,\eps_{0})$ and $t\in [0,T_{\eps}^{*})$.
\end{theorem}
\begin{proof}
Since the function $\{t\mapsto\Omega_\eps(t)\}$ defined in formula~\eqref{dualcontrol} is continuous, for any fixed $T_0>0$
and $\eps_0,\sigma_0>0$, we can define the time (recall here that $\Omega(0)={\mathcal O}(\eps^2)$ as $\eps\to 0$)
\begin{equation}
\label{Tepsdef}
T^*_\eps:=\sup\big\{t\in [0,T_0]:\,\Omega_\eps(s)\leq\sigma_0,\,\,\text{for all $s\in (0,t)$}\big\}>0,
\end{equation}
for all $\eps\in(0,\eps_0)$. Therefore, by choosing the numbers $\sigma_0$ and $\eps_0$ sufficiently small,
by virtue of Lemma~\ref{firstKey}, we conclude that 
\begin{equation*}
0\leq \E(\psi_{\vep}(t)) -\E(r) \leq C\Omega_\eps(t)+{\mathcal O}(\eps^2)\leq{\mathcal A},\quad
\text{for all $\eps\in(0,\eps_0)$ and $t\in [0,T^*_\eps)$}.
\end{equation*}
Since $\|\psi_\eps(t)\|_{L^2}=\|r\|_{L^2}$, we are in the right position to exploit
the stability property of ground states (Proposition~\ref{stabilita}).\ Hence,
there exist two families of uniformly bounded functions $\hat\theta_{\vep}:\R^+\to[0,2\pi)$ 
and $\hat y_{\vep}:\R^+\to\R^N$ such that
$$
\Big\|e^{-\frac{{\rm i}}{\vep}\xi(t)\cdot [\vep x+x(t)]} e^{-\iu A(x(t))\cdot x}
\,\phi_{\vep}(\vep x+x(t),t)-e^{\iu\hat\theta_\eps(t)} r\big(x+\hat y_{\vep}(t)\big)\Big\|_{H^1}^2\leq C\Omega_\eps(t)+{\mathcal O}(\eps^2),
$$
for all $\eps\in(0,\eps_0)$ and any $t\in [0,T^*_\eps)$. In turn, by rescaling and setting
$\theta_\eps(t):=\eps\hat\theta_\eps(t)$ and $y_\eps(t):=x(t)-\eps\hat y_\eps(t)$, we get
$$
\Big\|e^{-\frac{\iu}{\vep}\xi(t)\cdot x-\frac{\iu}{\eps} A(x(t))\cdot (x-x(t))}
\phi_\eps(x,t)-e^{{\frac{\iu}{\eps}}\theta_\eps(t)}r\Big(\frac{x-y_\vep(t)}{\vep}\Big)\Big\|_{\H_{\vep}}^2\leq C\Omega_\eps(t)+{\mathcal O}(\eps^2),
$$
for all $\eps\in(0,\eps_0)$ and $t\in [0,T^*_\eps)$, namely inequality~\eqref{quasifine}, concluding the proof.
\end{proof}

\medskip

\section{Mass and momentum identities}
\label{massmom-id}

In the following lemma we obtain two important identities satisfied by the equation. 
Only after completion of the present paper, that the author discovered the second identity 
was independently obtained in~\cite{selvit}. For the sake of self-containedness 
we include our proof, which uses the first identity and it is shorter.

\begin{lemma}\lbl{ideder}
Let $\phi_\eps$ be the solution to problem~\eqref{probMF} corresponding to the initial data~\eqref{initialD}.
Then we have the identity 
\begin{equation}
\label{derfi}
\frac1{\vep^{N}}\frac{\partial |\phi_{\vep}|^{2}}{\partial t}(x,t)=-{\rm div}_{x}\, p_{\vep}^A(x,t),\quad
x\in\R^N,\, t\in\R^+.
\end{equation}
Moreover, for all $t\in\R^+$, we have the identity 
\begin{equation}
\label{momenttid}
\int \frac{\partial p^A_{\vep}}{\partial t}(x,t)dx=
-\int p_{\vep}^A(x,t)\times B(x)dx -\int\nabla V(x)\frac{|\phi_{\vep}(x,t)|^{2}}{\eps^N}dx,
\end{equation}
where $B=\nabla\times A$ is the magnetic field associated with the vector potential $A$.
\end{lemma}

\begin{remark}
The momentum identity~\eqref{momenttid}, which plays an important r\^ole in our asymptotic analysis, can be thought 
as an extension of the so called {\em Ehrenfest's theorem} in presence of a magnetic field $B$.
\end{remark}

\begin{remark}
It follows from the momentum identity~\eqref{momenttid} that for the nonlinear Schr\"odinger equation with no magnetic 
field ($\nabla\times A=0$ in $\R^N$) and with a constant electric potential ($\nabla V=0$ in $\R^N$) the momentum 
$t\mapsto\int p^A_{\vep}(x,t)dx$ is a constant of motion.
\end{remark}

\begin{remark}
Concerning the addenda in the right-end side of~\eqref{momenttid},
in the semiclassical regime, by the upcoming Lemma~\ref{tec1}, as $\eps\to 0$,
$$
\int p_{\vep}^A(x,t)\times B(x)dx
+\int\nabla V(x)\frac{|\phi_{\vep}(x,t)|^{2}}{\eps^N}dx\sim m\xi(t)\times B(y_\eps(t))+m\nabla V(y_\eps(t)).
$$
We will show that $y_\eps(t)$ remains close to $x(t)$, for $\eps$ small (cf.\ Lemma~\ref{y}). Hence, from the right-hand side
of~\eqref{momenttid} the Newton equation~\eqref{DriveS} naturally emerges, ruling the dynamics of a particle
subjected to an electric force $F_e=-\nabla V(x(t))$ and to a magnetic force
$F_b=-v(t)\times B(x(t))$, being $v=\dot x$ the  velocity.
\end{remark}

\begin{proof}
By the exponential decay of $r(x),\partial_i r(x)$ and $\partial^2_{ij}r(x)$ given by~\eqref{expdecay} and the 
fact that $\|A\|_{C^1}<\infty$, the initial data~\eqref{initialD} belongs to $H^2_{A,\eps}$.
Hence, by the regularity (see Proposition~\ref{wellP}), it follows that $\phi_\eps(t)$ belongs
to  $H^1(\R^N;\C)\cap H^2_{A,\eps}$ for all $t>0$. By the standard  Calder\'on-Zygmund inequality
$\|\partial^2_{ij}\phi_\eps(t)\|_{L^2}\leq C\|\Delta\phi_\eps(t)\|_{L^2}$ for all $t$ (see e.g.~\cite[Corollary 9.10]{GT})
and since, again, $\|A\|_{C^1}<\infty$, for any $i,j=1,\dots,N$ we get
\begin{align*}
\eps^2\|\partial^2_{ij}\phi_\eps(t)\|_{L^2} &\leq C\|\eps^2\Delta\phi_\eps(t)\|_{L^2}
\leq C \|\big(\frac{\eps}{\iu}\nabla-A(x)\big)^2\phi_\eps(t)\|_{L^2}
+C\|A(x)\cdot\nabla \phi_\eps(t)\|_{L^2}   \\
& +C\||A(x)|^2\phi_\eps(t)\|_{L^2}+C\|{\rm div}_xA(x)\phi_\eps(t)\|_{L^2} \\
&\leq  C \|\phi_\eps(t)\|_{H^2_{A,\eps}}+C\|\phi_\eps(t)\|_{H^1} <\infty,
\end{align*}
for all $t>0$. Hence $\phi_\eps(t)\in H^2(\R^N;\C)$, for all $t>0$.
Set, for $j=1,\dots,N$,
$$
(p_{\vep}^A)_j(x,t)=\frac{1}{\eps^N}\im\big(\bar\phi_\eps(x,t)(\eps \partial_j\phi_\eps(x,t)-\iu A_j(x)\phi_\eps(x,t)\big).
$$
To prove identity~\eqref{derfi} notice that, on one hand, we have
\begin{align*}
-{\rm div}_{x}\,p_{\vep}^A(x,t)& =-\sum_{j=1}^N\partial_j(p_{\vep}^A)_j(x,t)  \\
& =-\sum_{j=1}^N\frac{1}{\eps^N}\im\big(\partial_j\bar\phi_\eps(x,t)(\eps \partial_j\phi_\eps(x,t)-\iu A_j(x)\phi_\eps(x,t)\big)   \\
& -\sum_{j=1}^N\frac{1}{\eps^N}\im\big(\bar\phi_\eps(x,t)(\eps \partial_{jj}^2\phi_\eps(x,t)-\iu \partial_jA_j(x)\phi_\eps(x,t)-\iu A_j(x)\partial_j\phi_\eps(x,t)\big) \\
& =\frac{2}{\eps^N}A(x)\cdot\re\big(\nabla\bar\phi_\eps(x,t)\phi_\eps(x,t)\big)   \\
& -\frac{1}{\eps^{N-1}}\im\big(\bar\phi_\eps(x,t)\Delta\phi_\eps(x,t))+\frac{1}{\eps^N}{\rm div}_x\,A(x)|\phi_\eps(x,t)|^2.
\end{align*}
On the other hand, it follows
\begin{align*}
\frac1{\vep^{N}}\frac{\partial |\phi_{\vep}|^{2}}{\partial t}(x,t)& =\frac2{\vep^{N+1}}\im\big(\bar\phi_{\vep}(x,t)\big[\frac{1}{2}\big(\frac{\eps}{\iu}\nabla-A(x)\big)^2\phi_\eps(x,t)
+V(x)\phi_\eps(x,t)-|\phi_\eps(x,t)|^{2p}\phi_\eps(x,t)\big]\big) \\
& =\frac1{\vep^{N+1}}\im\big(\bar\phi_{\vep}(x,t)\big(\frac{\eps}{\iu}\nabla-A(x)\big)^2\phi_\eps(x,t)\big) \\
& =-\frac1{\vep^{N-1}}\im\big(\bar\phi_{\vep}(x,t)\Delta\phi_\eps(x,t))
+\frac{2}{\vep^{N}}A(x)\cdot\re\big(\phi_{\vep}(x,t)\nabla\bar\phi_\eps(x,t)\big) \\
&+\frac{1}{\eps^N}{\rm div}_x\,A(x)|\phi_\eps(x,t)|^2.
\end{align*}
Now, concerning second identity, \eqref{momenttid}, for any $j=1,\dots,N$, it holds
\begin{align*}
\frac{\partial (p_{\vep}^A)_j}{\partial t} & =\vep^{1-N}\im(\partial_t\overline{\phi}_{\vep}\partial_j\phi_{\vep})
+\vep^{1-N}\im(\overline{\phi}_{\vep}\partial_j(\partial_t\phi_{\vep}))-\frac{1}{\eps^N}A_j(x)\frac{\partial |\phi_\eps|^2}{\partial t} \\
& =\vep^{1-N}\im(\partial_t\overline{\phi}_{\vep}\partial_j\phi_{\vep})
+\vep^{1-N}\im (\partial_j\big(\overline{\phi}_{\vep}\partial_t\phi_{\vep}\big))
- \vep^{1-N}\im (\partial_j\overline{\phi}_{\vep}\partial_t\phi_{\vep})-\frac{1}{\eps^N}A_j(x)\frac{\partial |\phi_\eps|^2}{\partial t} \\
\noalign{\vskip3pt}
& =2\vep^{1-N}\im(\partial_t\overline{\phi}_{\vep}\partial_j\phi_{\vep})
+\vep^{1-N}\im (\partial_j\big(\overline{\phi}_{\vep}\partial_t\phi_{\vep}\big))-\frac{1}{\eps^N}A_j(x)\frac{\partial |\phi_\eps|^2}{\partial t}.
\end{align*}
The second term integrates to zero. Moreover, 
taking into account identity~\eqref{derfi}, we get
\begin{align*}
-\int\frac{1}{\eps^N}A_j(x)\frac{\partial |\phi_\eps|^2}{\partial t}(x,t)dx &=\int A_j(x){\rm div}_{x}\, p_{\vep}^A(x,t)dx=
-\int \nabla A_j(x)\cdot p_{\vep}^A(x,t)dx \\
& =-\eps^{1-N} \int \sum_{i=1}^N \partial_i A_j(x)\im\big(\bar\phi_\eps(x,t)\partial_i\phi_\eps(x,t)\big)dx \\
& +\eps^{-N}\int \sum_{i=1}^N A_i(x)\partial_i A_j(x)|\phi_\eps(x,t)|^2dx.
\end{align*}
Concerning the first term in the formula for $\partial_t (p_{\vep}^A)_j$, conjugate the equation, 
multiply it by $2\eps^{-N}\partial_j\phi_\eps$ and take the imaginary part. It follows (summation on repeated $i$ indexes)
\begin{align*}
2\vep^{1-N} \im(\partial_{t}\overline{\phi_{\vep}}\partial_j\phi_{\vep}) &=-\vep^{2-N}\re(\Delta 
\overline{\phi_{\vep}}\partial_j\phi_{\vep})
+\eps^{-N}|A(x)|^2\re(\overline{\phi_{\vep}}\partial_j\phi_{\vep})  \\
\noalign{\vskip3pt}
& +\eps^{1-N}{\rm div}_xA(x)\im(\bar\phi_\eps\partial_j\phi_\eps)
+2\eps^{1-N}A(x)\cdot\im(\nabla \bar\phi_\eps\partial_j\phi_\eps) \\
\noalign{\vskip3pt}
& + 2\eps^{-N}V(x)\re(\overline{\phi_{\vep}}\partial_j\phi_{\vep}) 
-2\eps^{-N}|\phi_{\vep}|^{2p}\re(\overline{\phi_{\vep}}\partial_j\phi_{\vep})  \\
\noalign{\vskip3pt}
&=-\vep^{2-N}\re(\partial_i\big(\partial_i \overline{\phi_{\vep}}\partial_j\phi_{\vep}))+\vep^{2-N}\partial_j\Big(\frac{|\partial_i \phi_{\vep}|^{2}}{2}\Big) \\
\noalign{\vskip3pt}
& +\eps^{-N}|A(x)|^2\re(\overline{\phi_{\vep}}\partial_j\phi_{\vep})  
+\eps^{1-N}{\rm div}_xA(x)\im(\bar\phi_\eps\partial_j\phi_\eps) \\
\noalign{\vskip3pt}
& +2\eps^{1-N}A(x)\cdot\im(\nabla \bar\phi_\eps\partial_j\phi_\eps) + \eps^{-N}\partial_j\left(V(x) |\phi_{\vep}|^2\right)  \\
\noalign{\vskip3pt}
& -\eps^{-N}\partial_jV(x)|\phi_\vep|^2  -\eps^{-N}\frac{1}{p+1}\partial_j\big(|\phi_{\vep}|^{2p+2}\big).
\end{align*}
Notice that the following identity can be easily shown (recall that  $\phi_\eps(t)\in H^2$ for all $t$),
$$
\int {\rm div}_xA(x)\im(\bar\phi_\eps\partial_j\phi_\eps)dx+
2\!\!\int A(x)\cdot\im(\nabla \bar\phi_\eps\partial_j\phi_\eps)dx=
\int \sum_{i=1}^N \partial_j A_i(x)\im\big(\bar\phi_\eps\partial_i\phi_\eps\big)dx.
$$
Then, recalling that ${\mathbb H}^B=(\partial_j A_i-\partial_iA_j)_{ij}$ and that
the flux of $\phi_\eps$ is in $H^2$, we infer that
\begin{align*}
\int\frac{\partial (p_{\vep}^A)_j}{\partial t} & =-\eps^{-N}\int \sum_{i=1}^N A_i(x)\big (\partial_j A_i(x)-\partial_i A_j(x)) |\phi_{\vep}|^2dx  \\
\noalign{\vskip3pt}
& +\eps^{1-N}\int \sum_{i=1}^N \partial_jA_i(x)\cdot\im(\bar\phi_\eps\partial_i\phi_\eps)dx -\eps^{-N}\int \partial_jV(x)|\phi_\vep|^2dx  \\
& -\eps^{1-N} \int  \sum_{i=1}^N \partial_i A_j(x)\cdot\im\big(\bar\phi_\eps\partial_i\phi_\eps\big)dx   \\
& =-\eps^{-N}\int \sum_{i=1}^N A_i(x)\big (\partial_j A_i(x)-\partial_i A_j(x)) |\phi_{\vep}|^2dx  \\
\noalign{\vskip3pt}
& +\eps^{1-N}\int \sum_{i=1}^N (\partial_j A_i(x)-\partial_iA_j(x))\cdot\im(\bar\phi_\eps\partial_i\phi_\eps)dx -\eps^{-N}\int \partial_jV(x)|\phi_\vep|^2dx  \\
& =\int ({\mathbb H}^Bp^A_\eps(x,t))_jdx -
\eps^{-N}\int \partial_jV(x)|\phi_\vep|^2dx .
\end{align*}
Taking into account the formal identification of  the notation $-p_{\vep}^A(x,t)\times B(x)$ with
the matrix operation ${\mathbb H}^Bp^A_\eps(x,t)$, we obtain the assertion. 
To see this in the three dimensional case, recalling that 
$$
(B_1,B_2,B_3)=\nabla\times A=(\partial_2A_3-\partial_3A_2,\partial_3A_1-\partial_1A_3,\partial_1A_2-\partial_2A_1),
$$
we obtain the skew-symmetric matrix
\begin{equation*}
{\mathbb H}^B(x)=
\begin{bmatrix}
	0 & \partial_2 A_1-\partial_1A_2 & \partial_3 A_1-\partial_1A_3 \\
\partial_1 A_2-\partial_2A_1 & 0 & \partial_3 A_2-\partial_2A_3  \\
\partial_1 A_3-\partial_3A_1 & \partial_2 A_3-\partial_3A_2 & 0	
\end{bmatrix}
=
\begin{bmatrix}
	0 & -B_3   & B_2 \\
B_3   &  0     & -B_1  \\
-B_2  &  B_1   & 0	
\end{bmatrix}.
\end{equation*}
Then, setting $p_\eps^i=(p_{\vep}^A)_i$, it follows that
\begin{equation*}
{\mathbb H}^B(x)p_{\vep}^A(x,t)=
\begin{bmatrix}
	0 & -B_3   & B_2 \\
B_3   &  0     & -B_1  \\
-B_2  &  B_1   & 0	
\end{bmatrix}
\begin{bmatrix}
	p_{\vep}^1 \\
	p_{\vep}^2  \\
	p_{\vep}^3
\end{bmatrix}
=
\begin{bmatrix}
	p_\eps^3B_2-p_\eps^2B_3 \\
p_\eps^1B_3-p_\eps^3B_1 \\
p_\eps^2B_1-p_\eps^1B_2 
\end{bmatrix}
= -p_{\vep}^A(x,t)\times B(x).
\end{equation*}
The proof is now concluded.
\end{proof}

\medskip

\section{Mass and momentum estimates}
\label{massmomest}

First, we have the following control on the mass and momentum.

\begin{lemma}
	\label{tec1}
Let $\eps_{0}>0$, $T_{\eps}^{*}>0$ and $y_\eps(t)$ be as in Theorem~\ref{chiave}.
Then there exists a positive constant $C$ such that 
\begin{equation*}
\Big\|\frac{|\phi_{\vep}(x,t)|^2}{\vep^{N}}dx-m\delta_{y_{\vep}(t)}\Big\|_{(C^{2})^{*}}
+
\Big\|p_{\vep}^{A(x(t))}(x,t)dx-m\xi(t)\delta_{y_{\vep}(t)}\Big\|_{(C^{2})^{*}}\leq C\Omega_\eps(t)+{\mathcal O}(\eps^2),
\end{equation*}
for every $t\in[0,T_{\vep}^{*})$ and $\eps\in(0,\vep_{0})$.
\end{lemma}

\begin{proof}
For any $v\in H^{1}(\R^N;\C)$, we have $|\nabla |v||^{2}=|\nabla v|^{2}-\frac{|\im(\bar v\nabla v)|^{2}}{|v|^{2}}$.
Then, if $\psi_\eps(x,t)$ is the function introduced in formula~\eqref{veps}, by Lemma~\ref{firstKey} it follows that
\begin{equation*}
0\leq \E(|\psi_{\vep}|) -\E(r) +
\frac{1}{2}\int\frac{|\im(\bar \psi_{\vep}\nabla \psi_{\vep})|^{2}}{|\psi_{\vep}|^{2}}dx\leq C\Omega_\eps(t)+{\mathcal O}(\eps^2),
\end{equation*}
for every $t\in[0,T_{\vep}^{*})$ and $\eps\in(0,\vep_{0})$.
Moreover, as $\||\psi_{\vep}|\|_{L^2}=\|r\|_{L^2}$, by~\eqref{variatcaract-r} we have
\begin{equation}\label{dis1}
	\int\frac{|\im(\bar \psi_{\vep}\nabla \psi_{\vep})|^{2}}{|\psi_{\vep}|^{2}}dx\leq C\Omega_\eps(t)+{\mathcal O}(\eps^2),
\end{equation}
for every $t\in[0,T_{\vep}^{*})$ and $\eps\in(0,\vep_{0})$. Now, by the definition of $\psi_\eps$ (cf.~\eqref{veps}), we get
\begin{align*}
&\frac{|\im(\bar \psi_{\vep}\nabla \psi_{\vep})|^{2}}{|\psi_{\vep}|^{2}} \\
&=\frac{\left|\im(\bar\phi_{\vep}(\vep x+x(t),t)\eps\nabla\phi_{\vep}(\vep x+x(t),t))-(\xi(t)+A(x(t))|\phi_{\vep}(\vep x+x(t),t)|^{2}\right|^{2}}{|\phi_{\vep}(\vep x+x(t),t)|^{2}}
\\
&=\frac{\big|\eps^Np^{A(x(t))}_\eps(\eps x+x(t),t)-\xi(t)|\phi_{\vep}(\vep x+x(t),t)|^{2}\big|^{2}}{|\phi_{\vep}(\vep x+x(t),t)|^{2}}
\\
&=\eps^{2N}\frac{\big|p^{A(x(t))}_\eps(\eps x+x(t),t)\big|^{2}}{|\phi_{\vep}(\vep x+x(t),t)|^{2}}
+|\xi(t)|^{2}|\phi_{\vep}(\vep x+x(t),t)|^{2}  \\
\noalign{\vskip4pt}
&-2\eps^N\xi(t)\cdot p^{A(x(t))}_\eps(\eps x+x(t),t).
\end{align*}
Hence, by a change of variable, we reach
\begin{equation}
	\label{equa1}
\int \frac{|\im(\bar \psi_{\vep}\nabla \psi_{\vep})|^{2}}{|\psi_{\vep}|^{2}}dx
=\vep^{N}\int \frac{\big|p^{A(x(t))}_\eps(x,t)\big|^{2}}{|\phi_{\vep}(x,t)|^{2}}dx
+m|\xi(t)|^2 -2\xi(t)\cdot \int p_{\vep}^{A(x(t))}(x,t)dx.
\end{equation}
Notice that by simple computations, by combining~\eqref{dis1} and~\eqref{equa1}, it holds
\begin{align}\label{dis2}
 \int\Big|\vep^{N/2}\frac{p^{A(x(t))}_\eps(x,t)}{|\phi_{\vep}(x,t)|}
-\frac{ \int p^{A(x(t))}_\eps(x,t)dx}{m}
\frac{|\phi_{\vep}(x,t)|}{\vep^{N/2}}\Big|^{2}
&+m\Big|\xi(t)
-\frac{ \int p^{A(x(t))}_\eps(x,t)dx}{m}\Big|^{2}
\\
\nonumber \leq C\Omega_\eps(t)+{\mathcal O}(\eps^2),  &
\end{align}
for every $t\in[0,T_{\vep}^{*})$ and $\eps\in(0,\vep_{0})$.
To prove the assertion, we estimate  $\rho_{\vep}(t)$, where 
\beq\lbl{somma}
\rho_{\vep}(t):= \Big|\int \psi(x)\frac{|\phi_{\vep}(x,t)|^{2}}{\eps^N}dx-m\psi(y_{\vep})\Big|+
\Big|\int p^{A(x(t))}_\eps(x,t)\psi(x)-m\xi(t)\psi(y_{\vep})\Big|
\eeq
for every function $\psi$ of class $C^{2}$ such that $\|\psi\|_{C^{2}}\leq 1$. Taking into account that, by the definition of 
$\Omega_\eps(t)$ (cf. formula~\eqref{defOmega}) and $\|A-A(x(t))\|_{C^3}\leq C$, we have
\begin{align*}
\Big|\int p^{A(x(t))}_\eps(x,t)dx-m\xi(t)\Big|  & \leq 
\Big|\int p^{A}_\eps(x,t)dx-m\xi(t)\Big|+\Big|\int (p^{A}_\eps(x,t)-p^{A(x(t))}_\eps(x,t))dx\Big| \\
& \leq C\Omega_\eps(t)+\Big|\int (A(x)-A(x(t))) \frac{|\phi_\eps(x,t)|^2}{\eps^N}dx\Big| \\
& =  C\Omega_\eps(t)+\Big|\int (A(x)-A(x(t))) \Big(\frac{|\phi_\eps(x,t)|^2}{\eps^N}-m\delta_{x(t)}\Big)dx\Big| \\
& \leq C\Omega_\eps(t),
\end{align*}
we can conclude that
\begin{align*}
&\Big|\int p^{A(x(t))}_\eps(x,t)\psi(x)dx-m\xi(t)\psi(y_{\vep}(t))\Big|\leq \\ 
&\leq \Big|\int p^{A(x(t))}_\eps(x,t)[\psi(x)-\psi(y_{\vep}(t))]dx\Big|+
|\psi(y_{\vep}(t))|\Big|\int p^{A(x(t))}_\eps(x,t)dx-m\xi(t)\Big|
\\
&\leq \Big|\int p^{A(x(t))}_\eps(x,t)[\psi(x)-\psi(y_{\vep}(t))]dx\Big|+C\Omega_\eps(t)
\\
&\leq \frac1{m}\Big|\int p^{A(x(t))}_\eps(x,t)dx\Big|\Big|\int\frac{\psi(x)|\phi_{\vep}
(x,t)|^{2}}{\vep^{N}}dx-m\psi(y_{\vep}(t))\Big|
\\
&\;\,+\Big|\int \psi(x)\Big[p^{A(x(t))}_\eps(x,t)-\frac1{m}\Big(\int p^{A(x(t))}_\eps(x,t)dx\Big) 
\frac{|\phi_{\vep}(x,t)|^{2}}{\vep^{N}}\Big]dx\Big|+C\Omega_\eps(t),
\end{align*}
for all $\eps\in(0,\eps_{0})$ and $t\in [0,T_{\eps}^{*})$. 
Since $\int p^{A(x(t))}_\eps(x,t)dx$ is bounded (see Lemma~\ref{boundedMoment}) and  
$$
\int \Big[p^{A(x(t))}_\eps(x,t)-\frac1{m}\big(\int p^{A(x(t))}_\eps(x,t)dx\big)
\frac{|\phi_{\vep}(x,t)|^2}{\vep^{N}}\Big]dx=0,
$$
setting $\hat \psi(x):=\psi(x)-\psi(y_{\vep}(t))$, it holds
\begin{align*}
\rho_{\vep}(t) &\leq \int |\hat\psi(x)|\frac{|\phi_{\vep}(x,t)|^{2}}{\eps^N}dx+
C\int |\hat\psi(x)|\frac{|\phi_{\vep}(x,t)|^{2}}{\eps^N}dx \\
&+\int|\hat\psi(x)|\Big|p^{A(x(t))}_\eps(x,t)
-\frac1{m}\Big(\int p^{A(x(t))}_\eps(x,t)dx\Big)\frac{|\phi_{\vep}(x,t)|^{2}}{\vep^{N}}\Big|dx+C\Omega_\eps(t).
\end{align*}
From Young inequality and estimate~\eqref{dis2}, it follows 
\begin{align}
	\label{ineqqcrosscit}
 \rho_{\vep}(t) & \leq \int \big[ C|\hat\psi(x)| +\frac12|\hat\psi(x)|^{2}\big]\frac{|\phi_{\vep}(x,t)|^{2}}{\eps^N}dx  \\
\noalign{\vskip2pt}
&+ \frac12\int \Big|\vep^{N/2}\frac{p^{A(x(t))}_\eps(x,t)}{|\phi_{\vep}(x,t)|}
-\frac{1}{m}\Big(\int p^{A(x(t))}_\eps(x,t)dx\Big)\frac{|\phi_{\vep}(x,t)|}{\vep^{N/2}}\Big|^{2}+C\Omega_\eps(t)  \notag \\
\noalign{\vskip3pt}
&\leq  \int\big[ C|\hat\psi(x)| +\frac12|\hat\psi(x)|^{2}\big]\frac{|\phi_{\vep}(x,t)|^{2}}{\eps^N}dx+C\Omega_\eps(t)+{\mathcal O}(\eps^2). \notag
\end{align}
Via inequality $a^{2}\leq 2b^{2}+2(a-b)^{2}$ with $a=\eps^{-N/2}|\phi_{\vep}(x,t)|$ 
and $b=\eps^{-N/2}r((x-y_{\vep}(t))/\vep)$, 
\begin{align*}
\rho_{\vep}(t)  &\leq \frac{C}{\vep^{N}}\int\big[|\hat\psi(x)| +|\hat \psi(x)|^{2}\big]
r^{2}\big(\frac{x-y_{\vep}(t)}{\vep}
\big)dx +\frac{C}{\vep^{N}}\int \Big||\phi_{\vep}(x,t)|-r\big(\frac{x-y_{\vep}(t)}{\vep}\big)\Big|^{2}dx \\
\noalign{\vskip3pt}
&+C\Omega_\eps(t)+{\mathcal O}(\eps^2)\leq \Omega_\eps(t)+{\mathcal O}(\eps^2),
\end{align*}
for all $\eps\in(0,\eps_{0})$ and $t\in [0,T_{\eps}^{*})$,
by Lemma~\ref{pote} (as $\hat\psi (y_{\vep}(t))=0$) and Theorem~\ref{chiave}.
\end{proof}
\vskip2pt

Next, we need to show that the distance between the points $y_{\eps}(t)$ found out in the proof of 
Theorem~\ref{chiave} and the trajectory $x(t)$ is controlled by $\Omega_\eps(t)$, as $\eps$ goes to zero.

\begin{remark}
We stress that in the proof of the next Lemma we will choose the value of $T_0$ that was introduced 
in formula~\eqref{Tepsdef} inside the definition of $T^*_\eps$.
\end{remark}

\begin{lemma}\lbl{y}
Let $y_\eps(t)$ be as in Theorem~\ref{chiave}.
There exist positive constants $\eps_{0}$, $\sigma_{0}$ and $T_0$, namely the values introduced in~\eqref{Tepsdef}
in the definition of $T^*_\eps$ such that, for some positive constant $C$,
$$
|x(t)-y_{\vep}(t)|\leq C\Omega_\eps(t)+{\mathcal O}(\eps^2),
$$
for all $t\in [0,T_{\vep}^{*})$ and $\eps\in(0,\eps_{0})$.
\end{lemma}
\begin{proof}
We first show that there exists a time $T_{0}$ such that $|y_{\vep}(t)|<\rho$, for 
every $t\in [0,T_{\vep}^{*})$ with $T_{\vep}^{*}\leq T_{0}$, where $\rho$ is the positive 
constant introduced in formula~\eqref{defa}.  Let us first prove that $\|\delta_{y_{\vep}(t_{2})}-\delta_{y_{\vep}(t_{1})}\|_{C^{2*}}<\rho$
for all $t_{1},\,t_{2}\in [0,T_{\vep}^{*})$. Let $\varphi\in C^2(\R^N)$ be such that $\|\varphi\|_{C^2}\leq 1$. Hence, 
taking into account Lemma~\ref{boundedMoment} and identity~\eqref{derfi}, we get
		\begin{align*}
		\int\Big(\frac{|\phi_\eps(x,t_2)|^2}{\eps^N}-\frac{|\phi_\eps(x,t_1)|^2}{\eps^N}\Big)\varphi(x) dx 
		& =\int\int_{t_1}^{t_2}\frac1{\vep^{N}}\frac{\partial |\phi_{\vep}|^{2}}{\partial t}(x,t)\varphi(x) dtdx \\
		 =\int\int_{t_1}^{t_2}-\varphi(x){\rm div}_{x}\, p_{\vep}^A(x,t)dtdx 
&		 =\int_{t_1}^{t_2}\int\nabla\varphi(x)\cdot p_{\vep}^A(x,t)dxdt \\
		\leq \|\nabla \varphi\|_{L^\infty}\int_{t_1}^{t_2}dt\int |p_{\vep}^A(x,t)|dx &\leq C\|\varphi\|_{C^2}|t_2-t_1|\leq C|t_2-t_1|.
	\end{align*}
Hence, for all $t_{1},\,t_{2}\in [0,T_{\vep}^{*})$, it holds
$$
\Big\| \frac{|\phi_\eps(x,t_2)|^2}{\eps^N}dx-\frac{|\phi_\eps(x,t_1)|^2}{\eps^N}dx \Big\|_{{C^{2*}}}\leq C|t_2-t_1|.
$$
In view of Lemma~\ref{tec1}, the following inequality holds, 
$$
m\|\delta_{y_{\vep}(t_{2})}-\delta_{y_{\vep}(t_{1})}\|_{C^{2*}}\leq CT_{0}+C\Omega_\eps(t)+{\mathcal O}(\eps^2)\leq C(T_0+\sigma_0)+{\mathcal O}(\eps^2).
$$
Here we choose the value of $T_{0}$ and then of $\sigma_{0},\eps_0$ so small that
$$C(T_0+\sigma_0)+{\mathcal O}(\eps^2)<\min\{mK_0,mK_0K_1\},$$ being $K_0$ and $K_1$ the constants
introduced in Lemma~\ref{propdelta}. Hence, $|y_{\vep}(t_{2})-y_{\vep}(t_{1})|< K_{0}$
for all $t_{1},\,t_{2}\in [0,T_{\vep}^{*})$, and since $y_{\vep}(0)=x_0$,
we obtain the desired assertion. We can now conclude the proof of this Lemma.
The properties of the function $\chi$ imply
\begin{equation*}
|x(t)-y_{\vep}(t)|\leq \frac1{m}|\gamma_{\vep}(t)|+\frac1{m}
\Big|\int x\chi(x)\frac{|\phi_{\vep}(x,t)|^{2}}{\eps^N}dx-my_{\vep}(t)\Big|.
\end{equation*}
In light of the first step of the proof, we have $\chi(y_{\vep}(t))=1$ for all $t\in [0,T_{\vep}^{*})$ 
and $\eps\in(0,\eps_{0})$, so that exploiting again Lemma~\ref{tec1}, we conclude that
\begin{equation*}
|x(t)-y_{\vep}(t)|\leq C\Omega_\eps(t)+C\|x\chi\|_{C^{2}}\Big\|
\frac{|\phi_{\vep}(x,t)|^{2}}{\vep^{N}}\,dx-m\delta_{y_{\vep}(t)}\Big\|_{C^{2*}}\leq C\Omega_\eps(t)+{\mathcal O}(\eps^2),
\end{equation*}
which yields the assertion.
\end{proof}

Finally, we get a strengthened version of Lemma~\ref{tec1}.

\begin{lemma}
\label{tec1+}
Let $\eps_{0}>0$ and $T_{\eps}^{*}>0$ be as in Theorem~\ref{chiave}.
Then there exists a positive constant $C$ such that 
\begin{equation*}
\Big\|\frac{|\phi_{\vep}(x,t)|^2}{\vep^{N}}dx-m\delta_{x(t)}\Big\|_{C^{2*}}
+
\big\|p_{\vep}^{A}(x,t)dx-m\xi(t)\delta_{x(t)}\big\|_{C^{2*}}\leq C\Omega_\eps(t)+{\mathcal O}(\eps^2),
\end{equation*}
for every $t\in[0,T_{\vep}^{*})$ and $\eps\in(0,\vep_{0})$. In particular, if $\|A\|_{C^2}$ is sufficiently small, we have
\begin{equation}
	\label{mom-fin-ineqq}
\Big\|\frac{|\phi_{\vep}(x,t)|^2}{\vep^{N}}dx-m\delta_{x(t)}\Big\|_{C^{2*}}
+
\big\|p_{\vep}^{A}(x,t)dx-m\xi(t)\delta_{x(t)}\big\|_{C^{2*}}\leq C\hat\Omega_\eps(t)+{\mathcal O}(\eps^2),
\end{equation}
for every $t\in[0,T_{\vep}^{*})$ and $\eps\in(0,\vep_{0})$.
\end{lemma}
\begin{proof}
Notice that, taking into account Lemma~\ref{tec1}, Lemma~\ref{propdelta} and Lemma~\ref{y}, we get
\begin{align*}
\Big\|\frac{|\phi_{\vep}(x,t)|^2}{\vep^{N}}dx-m\delta_{x(t)}\Big\|_{C^{2*}}& \leq
\Big\|\frac{|\phi_{\vep}(x,t)|^2}{\vep^{N}}dx-m\delta_{y_\eps(t)}\Big\|_{C^{2*}} \\
\noalign{\vskip3pt}
&+m\big\|\delta_{y_\eps(t)}-\delta_{x(t)}\big\|_{C^{2*}}   
\leq C\Omega_\eps(t)+{\mathcal O}(\eps^2),
\end{align*}
for every $t\in[0,T_{\vep}^{*})$ and $\eps\in(0,\vep_{0})$. In turn, we also get
\begin{align*}
& \big\|p_{\vep}^A(x,t)dx-m\xi(t)\delta_{x(t)}\big\|_{C^{2*}} \leq \big\|p_{\vep}^A(x,t)dx-p^{A(x(t))}_\eps(x,t)dx\big\|_{C^{2*}} \\
\noalign{\vskip5pt}
&+\big\|p^{A(x(t))}_\eps(x,t)dx-m\xi(t)\delta_{y_{\vep}(t)}\big\|_{C^{2*}} +\big\|m\xi(t)\delta_{y_{\vep}(t)}-m\xi(t)\delta_{x(t)}\big\|_{C^{2*}}   \\
\noalign{\vskip2pt}
&\leq  \sup_{\|\varphi\|_{C^2}\leq 1}\Big|\int [A(x)-A(x(t))]\varphi(x)\frac{|\phi_\eps(x,t)|^2}{\eps^N}dx\Big|+C\Omega_\eps(t)+{\mathcal O}(\eps^2) \\
& =\sup_{\|\varphi\|_{C^2}\leq 1}\Big|\int [A(x)-A(x(t))]\varphi(x)\Big[\frac{|\phi_{\vep}(x,t)|^2}{\vep^{N}}dx-m\delta_{x(t)}\Big]dx\Big|
+C\Omega_\eps(t)+{\mathcal O}(\eps^2) \\
& \leq\sup_{\|\varphi\|_{C^2}\leq 1}\|(A(x)-A(x(t)))\varphi(x)\|_{C^2}\Big\|\frac{|\phi_{\vep}(x,t)|^2}{\vep^{N}}dx-m\delta_{x(t)}\Big\|_{C^{2*}}
+C\Omega_\eps(t)+{\mathcal O}(\eps^2) \\
\noalign{\vskip2pt}
&\leq C\Big\|\frac{|\phi_{\vep}(x,t)|^2}{\vep^{N}}dx-m\delta_{x(t)}\Big\|_{C^{2*}}
+C\Omega_\eps(t)+{\mathcal O}(\eps^2)\leq  C\Omega_\eps(t)+{\mathcal O}(\eps^2),
\end{align*}
for every $t\in[0,T_{\vep}^{*})$ and $\eps\in(0,\vep_{0})$. This concludes the proof
of the first assertion. Taking into account the definitions of $\Omega_\eps(t)$ and $\rho_\eps^A(t)$,
inequality~\eqref{mom-fin-ineqq} is just a simple consequence.
\end{proof}

\section{Proof of the main result concluded}
\label{proofmain}

In this section we will conclude the proof of the main result.

\subsection{The error estimate}
\label{small-err}

We now show that the quantity $\Omega_\eps(t)$, introduced in~\eqref{defOmega}, can be made small
at the order ${\mathcal O}(\eps^2)$, uniformly on finite time intervals, as $\eps\to 0$.
\vskip4pt

\begin{lemma}
\label{smallstatem}
There exists a positive constant $C=C(T_0)$ such that
$\hat\Omega_\eps(t)\leq C(T_0)\eps^2$,
for all $\eps\in (0,\eps_0)$ and $t\in [0,T^*_\eps)$.
If in addition we assume that $\|A\|_{C^2}<\delta$ for some $\delta>0$ small, then
there exists a positive constant $C=C(T_0)$ such that
$\Omega_\eps(t)\leq C(T_0)\eps^2$,
for all $\eps\in (0,\eps_0)$ and $t\in [0,T^*_\eps)$.
\end{lemma}
\begin{proof}
Taking into account Lemma~\ref{tec1+}, via identity~\eqref{momenttid} of Lemma~\ref{ideder}, we obtain
	\begin{align*}
	 &\Big|\int \frac{d}{dt}\Pi^1_\eps(x,t)dx\Big|  
	=\Big|\int \frac{\partial p^A_\eps}{\partial t}(x,t)dx-m\dot \xi(t)\Big| \\
	&=\Big|\int p_{\vep}^A(x,t)\times B(x)dx+\int\nabla V(x)\frac{|\phi_{\vep}(x,t)|^{2}}{\eps^N}dx   \\
&	-m\nabla V(x(t))-m\xi(t)\times B(x(t))\Big|  \\
	&=\Big|\int p_{\vep}^A(x,t)\times B(x)dx+\int\nabla V(x)\frac{|\phi_{\vep}(x,t)|^{2}}{\eps^N}dx \\
&	-\int m\nabla V(x)\delta_{x(t)}dx-m\int \xi(t)\times B(x)\delta_{x(t)}dx\Big|   \\
	&\leq  \Big|\int \Big(p_{\vep}^A(x,t)-m\xi(t)\delta_{x(t)}\Big)\times B(x)dx\Big|  \\
	& +\Big|\int\nabla V(x)\Big(\frac{|\phi_{\vep}(x,t)|^{2}}{\eps^N}-m\delta_{x(t)}\Big)dx\Big| \\
	&\leq \|A\|_{C^3}\Big\|p_{\vep}^A(x,t)dx-m\xi(t)\delta_{x(t)}\Big\|_{C^{2*}} \\
	&+\|V\|_{C^3}\Big\|  \frac{|\phi_{\vep}(x,t)|^{2}}{\eps^N}dx-m\delta_{x(t)} \Big\|_{C^{2*}}  \\
	\noalign{\vskip3pt}
	&\leq C\hat\Omega_\eps(t)+{\mathcal O}(\eps^2),
	\end{align*}
	for all $\eps\in (0,\eps_0)$ and $t\in [0,T^*_\eps)$. Hence, recalling Lemma~\ref{controestimateee}, it follows that
	\begin{align}
\label{conlquasi}
	\left|\int \Pi^1_\eps(x,t)dx\right| & \leq\left|\int \Pi^1_\eps(x,0)dx\right|+\int_{0}^t \left|\int\frac{d}{dt}\Pi^1_\eps(x,\tau)dx\right|d\tau  \notag \\
	& \leq {\mathcal O}(\eps^2)+C\int_{0}^t \hat\Omega_\eps(\tau)d\tau.  
\end{align}
Let now $\varphi\in C^3(\R^N)$ with $\|\varphi\|_{C^3(\R^N)}\leq 1$. Then identity~\eqref{derfi} and Lemma~\ref{tec1+} yield
\begin{align*}
	 &\Big|\int \frac{d}{dt}\Pi^2_\eps(x,t)\varphi(x)dx\Big|  
	=\Big|\int \varphi\frac{\partial}{\partial t}\frac{|\phi_\eps(x,t)|^2}{\eps^N}dx-m\nabla\varphi(x(t))\cdot\xi(t)\Big| \\
	& =\Big|- \int \varphi(x)\,{\rm div}_{x}\, p_{\vep}^A(x,t)dx-m\nabla\varphi(x(t))\cdot\xi(t)\Big| \\
	& =\Big|\int \nabla \varphi(x)\cdot p_{\vep}^A(x,t)dx-\int m\nabla\varphi(x)\cdot\xi(t)\delta_{x(t)}dx\Big| \\
	& =\Big|\int \nabla \varphi(x)\cdot\Big (p_{\vep}^A(x,t)-m\xi(t)\delta_{x(t)}\Big)dx\Big| \\
	&\leq  \|\varphi\|_{C^3}\big\|p_{\vep}^A(x,t)dx-m\xi(t)\delta_{x(t)}\big\|_{C^{2*}}
	\leq C\hat\Omega_\eps(t)+{\mathcal O}(\eps^2),
	\end{align*}
for all $\eps\in (0,\eps_0)$ and $t\in [0,T^*_\eps)$. Hence, by Lemma~\ref{controestimateee}, it follows that
	\begin{align}
	\label{seconpezz}
	\sup_{\|\varphi\|_{C^3}\leq 1}\Big|\int \Pi^2_\eps(x,t)\varphi(x)dx\Big| & 
	\leq \sup_{\|\varphi\|_{C^3}\leq 1}\Big|\int \Pi^2_\eps(x,0)\varphi(x)dx\Big|  \\
	&+\sup_{\|\varphi\|_{C^3}\leq 1}
	\int_{0}^t \Big|\int\frac{d}{dt}\Pi^2_\eps(x,\tau)\varphi(x)dx\Big|d\tau  \notag \\
	& \leq {\mathcal O}(\eps^2)+C\int_{0}^t \hat\Omega_\eps(\tau)d\tau,  \notag
\end{align}
for all $\eps\in (0,\eps_0)$ and $t\in [0,T^*_\eps)$. Finally, again via identity~\eqref{derfi} and Lemma~\ref{tec1+},
\begin{align*}
\big|\dot\gamma_{\vep}(t)\big| &=\Big|m\xi(t)+\int x\chi(x){\rm div}_x\, p^A_\eps(x,t)dx\Big|  \\
&=\Big|m\xi(t)-\int \nabla(x\chi(x))\cdot p^A_\eps(x,t)dx\Big|  \\
&=\Big|\int \nabla(x\chi(x))m\xi(t)\delta_{x(t)}-\int \nabla(x\chi(x))\cdot p^A_\eps(x,t)dx\Big|  \\
\noalign{\vskip2pt}
&\leq \|\nabla(x\chi(x))\|_{C^2}\big\|p^A_\eps(x,t)dx-m\xi(t)\delta_{x(t)}\big\|_{C^{2*}} 
\leq C\hat\Omega_\eps(t)+{\mathcal O}(\eps^2),
\end{align*}
for all $\eps\in (0,\eps_0)$ and $t\in [0,T^*_\eps)$. This, recalling Lemma~\ref{controestimateee}, yields
\begin{equation}
\label{lastpezz}
|\gamma_\eps(t)|\leq {\mathcal O}(\eps^2)+C\int_0^t \hat\Omega_\eps(\tau)d\tau,
\end{equation}
for all $\eps\in (0,\eps_0)$ and $t\in [0,T^*_\eps)$.
By collecting inequalities~\eqref{conlquasi},~\eqref{seconpezz} and~\eqref{lastpezz}, we get 
$$
\hat\Omega_\eps(t)\leq {\mathcal O}(\eps^2)+C\int_{0}^t \hat\Omega_\eps(\tau)d\tau
$$
for all $\eps\in (0,\eps_0)$ and $t\in [0,T^*_\eps)$. Then, by Gronwall Lemma, we have 
$\hat\Omega_\eps(t)\leq C(T_0)\eps^2$, for all $\eps\in (0,\eps_0)$ and $t\in [0,T^*_\eps)$.
Finally, recalling the definitions of $\Omega_\eps(t)$ and $\rho_\eps^A(t)$
and exploiting again Lemma~\ref{tec1+}
concludes the proof.
\end{proof}

\vskip3pt
We are now ready to conclude the proof of Theorem~\ref{mainthBest}.
\vskip4pt
\noindent
Let $\delta>0$ be as in Lemma~\ref{smallstatem}. Let us prove the first part of Theorem~\ref{mainthBest}.

\noindent
We recall that the value of $T_0>0$ was fixed in the proof of Lemma~\ref{y} and it just depends on the data
of the problem, such as $V,A,m,N$. Moreover, by virtue of 
Lemma~\ref{smallstatem} and by the definition of $T^*_\eps$ (see the proof of Theorem~\ref{chiave}),
it follows that $T^*_\eps=T_0$ for all $\eps\in (0,\eps_0)$, up to further reducing the value of $\eps_0$. 
Hence $\Omega_\eps(t)\leq C(T_0)\eps^2$  for all $\eps\in (0,\eps_0)$ and $t\in [0,T_0]$.
Now, by Theorem~\ref{chiave} there exist two families of functions $\theta_{\vep}:\R^+\to[0,2\pi)$ and  
$y_{\vep}:\R^+\to\R^N$ such that
\begin{equation*}
\Big\|	\phi_\eps(\cdot,t)-e^{\frac{\iu}{\vep}(\xi(t)\cdot x+\theta_\eps(t)+A(x(t))
	\cdot (x-x(t))}r\Big(\frac{x-y_\vep(t)}{\vep}\Big)\Big\|_{\H_\eps}^2={\mathcal O}(\eps^2), 
\end{equation*}	
for all $t\in [0,T_0]$. On the other hand, by combining  Lemma~\ref{y} with Lemma~\ref{smallstatem},
it follows that $|x(t)-y_{\vep}(t)|\leq C\eps^2$, for all $t\in [0,T_0]$ and $\eps\in(0,\eps_{0})$. Then,
taking into account the exponential decay of $\nabla r$, we obtain 
$$
\Big\|r\Big(\frac{x-y_\vep(t)}{\vep}\Big)-r\Big(\frac{x-x(t)}{\vep}\Big)\Big\|^2_{\H_\eps}\leq C\frac{|x(t)-y_\eps(t)|^2}{\eps^2}={\mathcal O}(\eps^2),
$$
for all $t\in [0,T_0]$ and $\eps\in(0,\eps_{0})$. Therefore, Theorem~\ref{mainthBest} holds true on the time interval $[0,T_0]$.
Let us take $x(T_0)$ and $\xi(T_0)$ as new initial data in system~\eqref{DriveS} and the function
\begin{equation*}
\phi_0^{{\rm new}}(x):=r\Big(\frac{x-x(T_0)}{\eps}\Big)
e^{\frac{\iu}{\eps}[A(x(T_0))\cdot(x-x(T_0))+x\cdot\xi(T_0)]},
\end{equation*}
as a new initial data for problem~\eqref{probMF}. Whence, by the previous step of the proof,
the approximation result holds on the interval $[T_0,2T_0]$, and hence on an arbitrary finite 
time interval $[0,T]$, for $T>0$.
\vskip4pt
\noindent
In order to prove the second part of the statement of Theorem~\ref{mainthBest} one
can follow the argument of~\cite{selvit} (essentially relying on~\cite{bronski}).
Based upon the identity
$$
\Big|\frac{\nabla v}{\iu}-Av\Big|^2=\frac{|p^{A}(v)|^2}{|v|^2}+|\nabla |v||^2,
\qquad
p^A(v):=\im\big(\bar v(\nabla v-\iu A v)\big),
$$
the energy functional $E_\eps$ rewrites as
$$
E_\eps(t)=E_\eps^{{\rm pot}}(t)+E_\eps^{{\rm b}}(t)+E_\eps^{{\rm k}}(t),
$$
where we have set
\begin{align*}
E_\eps^{{\rm pot}}(t)&:=\frac{1}{\eps^N}\int V(x)|\phi_\eps(x,t)|^2dx, \qquad \,\,
E_\eps^{{\rm k}}(t) :=	\frac{\eps^N}{2}\int \frac{|p^{A}_\eps(x,t)|^2}{|\phi_\eps(x,t)|^2}dx,   \\
E_\eps^{{\rm b}}(t) &:=\frac{1}{2\eps^N}\int 	|\nabla |\phi_\eps|(x,t)|^2dx
-\frac{1}{p+1}\frac{1}{\eps^N}\int |\phi_\eps(x,t)|^{2p+2}dx .
\end{align*}
Then, following the steps of the proof of~\cite[Lemma 3.5]{selvit} (on the basis of the
quantitative estimate of the expansion of $E_\eps$ up to a error of ${\mathcal O}(\eps^2)$,
cf.\ Lemma~\ref{estUNO}), we get
\begin{equation}
	\label{estcruc0}
0\leq E_\eps^{{\rm b}}(|\phi_\eps|)-E_\eps^{{\rm b}}(r)\leq C\hat\Omega_\eps(t)+{\mathcal O}(\eps^2),
\end{equation}
\begin{equation*}
0\leq E_\eps^{{\rm k}}(t)-\frac{1}{2}\frac{\big|\int p^{A}_\eps(x,t) \big|^2}{m}
\leq C\hat\Omega_\eps(t)+{\mathcal O}(\eps^2).
\end{equation*}
In turn, the second inequality easily yields
\begin{equation}
	\label{estcruc1}
\int \Big |\eps^{N/2}\frac{p^{A}_\eps(x,t)}{|\phi_\eps(x,t)|}
-\frac{\big(\int p^{A}_\eps(x,t)\big)}{m} \frac{|\phi_\eps(x,t)|}{\eps^{N/2}}\Big|^2dx
\leq C\hat\Omega_\eps(t)+{\mathcal O}(\eps^2).
\end{equation}
Once inequalities~\eqref{estcruc0}-\eqref{estcruc1} holds true, the assertion
can be proved by arguing as before. In fact, inequality~\eqref{estcruc0} yields
\begin{equation*}
\||\phi_\eps|-r\big(\frac{\cdot -y_\eps(t)}{\eps}\big) \|_{\H_\eps}^{2}
\leq C\hat\Omega_\eps(t)+{\mathcal O}(\eps^2),
\end{equation*}
for some $y_\eps(t)\in\R^N$. 
Instead, inequality~\eqref{estcruc1} allows to prove
inequality~\eqref{mom-fin-ineqq} of Lemma~\ref{tec1+}.

\bigskip	
	
\section{Conclusions}
\label{conclusions}	
We have analyzed the soliton dynamics features of subcritical (with respect to global well-posedness)
nonlinear Schr\"odinger equations in the semiclassical
regime under the effects of an external electromagnetic field, showing that the solutions
concentrate along a smooth curve $x(t):\R^+\to\R^N$ which is a parameterization of a solution of the classical Newton equation
involving a conservative electric force $F_e=-\nabla V(x(t))$ as well as the contribution of the Lorenz 
force $F_b=-\dot x(t)\times B(x(t))$, being $B=\nabla\times A$ the magnetic field. 
The main results improves the results of~\cite{selvit}, a recent contribution that the author discovered after completion
of the paper. The technique is based upon the use
of quantum (mass and energy for the PDE~\eqref{probMF}) and classical (\eqref{Hamilt} for the ODE~\eqref{NewtNew}) conservation laws,
on the lines of an argument introduced in J.\ Bronski and R.\ Jerrard in 2000 in~\cite{bronski} making
no use of a linearization procedure for the equation. On the other hand, the presence of the magnetic field introduces
new difficulties that have to be handled.
Finally, we wish to stress that our results are consistent
with the current literature regarding the analysis of particular classes solutions, 
such as the standing waves.

\bigskip
\bigskip

\end{document}